\documentclass{article}     
\usepackage{color}
\usepackage{array}

\usepackage[normalem]{ulem} 
\usepackage[margin=1in]{geometry}

\usepackage{url}

\usepackage[ansinew]{inputenc}
\usepackage{fancyhdr}
\usepackage[active]{srcltx}
\usepackage{amsmath,amssymb}
\usepackage[pdftex]{graphicx}
\usepackage{authblk}

\newcommand{\shcdot}{\!\cdot\!}
\newcommand{\shddot}{\!:\!}

\newcommand{\Rbb}{\mathbb{R}}

\newcommand{\cA}{{\cal A}}
\newcommand{\cP}{{\cal P}}

\newcommand{\cS}{{\cal S}}

\newcommand{\ttrain}{\text{train}}
\newcommand{\tmin}{\text{min}}
\newcommand{\tmax}{\text{max}}

\newcommand{\aPOD}{^{\tiny\text{POD}}}

\newcommand{\tmid}{\text{mid}}
\newcommand{\tDL}{\text{DL}}
\newcommand{\tPOD}{\text{POD}}
\newcommand{\tFOM}{\text{FOM}}

\newcommand{\micr}{\,$\mu$m\,}

\newcommand{\veps}{\varepsilon}

\newcommand{\bfm}[1]{\mathbf{#1}}
\newcommand{\bsn}[1]{\boldsymbol{#1}}

\newcommand{\bfzero}{\bfm{0}}

\newcommand{\bfB}{\bfm{B}}

\newcommand{\bff}{\bfm{f}}

\newcommand{\bfu}{\bfm{u}}
\newcommand{\bfw}{\bfm{w}}
\newcommand{\bfx}{\bfm{x}}

\newcommand{\bfS}{\bfm{S}}
\newcommand{\bfX}{\bfm{X}}
\newcommand{\bfV}{\bfm{V}}

\newcommand{\bfP}{\bfm{P}}

\newcommand{\bfF}{\bfm{F}}
\newcommand{\bfG}{\bfm{G}}
\newcommand{\bfU}{\bfm{U}}
\newcommand{\tm}{\bfm{1}}    
\newcommand{\zerou}{\bfm{0}}
\newcommand{\bfM}{\bfm{M}}
\newcommand{\bfN}{\bfm{N}}
\newcommand{\bfK}{\bfm{K}}
\newcommand{\bfZ}{\bfm{Z}}

\newcommand{\bfC}{\bfm{C}}
\newcommand{\bfH}{\bfm{H}}

\newcommand{\bfE}{\bfm{E}}

\newcommand{\bfSigma}{\bfm{\Sigma}}

\newcommand{\bftheta}{\bsn{\theta}}
\newcommand{\bfphi}{\bsn{\phi}}

\newcommand{\bfmu}{\bsn{\mu}}

\newcommand{\dOmega}{{\text d}\Omega}

\begin{document}

\title{Deep learning-based reduced order models for the real-time simulation of the nonlinear dynamics of microstructures}

\author[1]{Stefania Fresca}
\author[2]{Giorgio Gobat}
\author[3]{Patrick Fedeli}
\author[2]{Attilio Frangi}
\author[1]{Andrea Manzoni}

\affil[1]{MOX - Department of Mathematics,Politecnico di Milano, P.za Leonardo da Vinci 32, 20133 Milano, Italy}

\affil[2]{Department of Civil and Environmental Engineering, Politecnico di Milano, P.za Leonardo da Vinci 32, 20133 Milano, Italy}

\affil[3]{Analog and MEMS Group, STMicroelectronics, Via Tolomeo - 20007 Cornaredo (MI),Italy}

%
\maketitle
\section*{Abstract} 
We propose a non-intrusive Deep Learning-based Reduced Order Model (DL-ROM) 
capable of capturing the complex dynamics of mechanical systems showing inertia and geometric nonlinearities. 
In the first phase, a limited number of high fidelity snapshots are used to
generate a POD-Galerkin ROM which is subsequently exploited to generate the data, covering the whole parameter range, used in the training phase of the
DL-ROM. 
A convolutional autoencoder is employed
to map the system response onto a low-dimensional representation 
and, in parallel, to model the reduced nonlinear trial manifold.
The system dynamics on the manifold is described by means of a deep feedforward neural network that is trained together with the autoencoder.
The strategy is benchmarked against high fidelity solutions on a clamped-clamped beam
and on a real micromirror with
softening response and multiplicity of solutions.  
By comparing the different computational costs, we discuss the impressive gain in performance 
and show that the DL-ROM truly represents a real-time tool which can be profitably
and efficiently employed in complex system-level simulation procedures for design and optimisation purposes.




\section{Introduction}
\label{sec:intro}

Nonlinear modelling in solid and structural mechanics has received an impressive boost in recent years thanks to the increasing availability of computational resources. In particular, the nonlinear dynamics of Micro-Electro-Mechanical-Systems (MEMS) has attracted great attention as nonlinearities-related effects play an important role at the microscale and they intrinsically arise during the regular functioning of MEMS. 
Although usually avoided and considered a nuisance, nonlinear signatures can be efficiently exploited paving the way to extraordinary performances 
for both existing devices and new nonlinear-ameliorated applications/working principles. Seminal works have already demonstrated the potential of devices that exploit nonlinear phenomena in microstructures. For instance, nonlinearities improve the frequency stability over temperature fluctuations of resonant accelerometers \cite{Shin18} and of resonators designed for clock applications \cite{Defoort16}; the combination of mechanical (hardening) and electrostatic (softening) nonlinearities also widens the linear response of MEMS sensors \cite{Agarwal06,Rivlin12}. 
Internal resonance, frequently observed in nano-micro mechanical resonators, has been employed to stabilize the oscillation frequency of self-sustaining resonators \cite{Chen17}, and to improve the performances of Coriolis gyroscopes \cite{Sarrafan19}. Parametric resonance breaks into the engineering of filters \cite{Rhoads10}, flow sensors \cite{Torteman19} and mass sensors \cite{Zhang05}, as well as in boosting the sensitivity of gyroscopes \cite{Nitzan15,Polunin17}. 
The design of MEMS devices based on the co-existence of multiple stable solutions has inspired gas sensors \cite{Bouchaala16}, shock sensors \cite{frangi2015}, accelerometers \cite{HaLevy20} and flow sensors \cite{Kessler16}. In order to properly master and exploit nonlinearities, improved modelling capabilities are essential.

In this paper, we focus on the large family of MEMS resonators (e.g., accelerometers, gyroscopes, magnetometers, micromirrors, time-keeping devices) which are primarily defined by the steady-state periodic response. Selected output quantities of interest like, e.g., the maximum midspan deflection of a beam, or the rotation amplitude of a micromirror, should be predicted as a function of the actuation intensity and frequency -- thus playing the role of input parameters -- and ultimately yielding the so-called Frequency Response Function (FRF).  Moreover, the actuation can be electrostatic, piezoelectric, or magnetic, according to the considered applications, hence introducing additional sources of nonlinearities. Tracing the FRF thus requires to solve the problem multiple times (ranging from $O(10^3)$ to $O(10^7)$, depending upon the application at hand) for different values of a set of input parameters. 

Even if numerical methods are emerging as a general solution to tackle the aforementioned tasks, the computational cost of high-fidelity approaches still remains a major issue, if not a true bottleneck. For instance, dedicated Harmonic Balance techniques or shooting procedures are overwhelmingly complex and time consuming \cite{Kerschen09,KerschenNNM}. On the other hand, because of the low dissipation in vacuum packaged devices, ultimately leading to long transients, time marching schemes are hardly computationally affordable, yielding to high-dimensional, nonlinear dynamical systems to be solved. The community of nonlinear dynamics at the microscale therefore lacks of approaches that are suitable for the new generation of MEMS, in specific scenarios where multi-queries and real-time performances are required for the sake of design, testing and optimization. In all these cases, relying on high-fidelity computational techniques is indeed unaffordable.
 
Among different options, reduced order models (ROMs) represent a key numerical tool in order to generate efficient, yet reliable, approximations to the solution of parametrized differential problems. In particular, projection-based ROMs reshape the original high-fidelity problem into a nonlinear, dynamical system featuring a much lower dimension, yet capable to capture the physical features of the problem at hand
\cite{mignolet,Besselink13,alfiobook}. However, projection-based ROMs  are intrusive techniques, possibly expensive to be assembled in the case of high order polynomial (or nonpolynomial) nonlinearities, and hurdles to scale over many parameters instances \cite{Manzoni17,Pagani18,Farhat20}.  

In this framework, extensions of the modal methods available for linear problems to the nonlinear regime are receiving increasing attention \cite{touze2006}. Nevertheless, the computation of the invariant manifolds of nonlinear normal modes starting from the normal form theory has been developed to maturity only for relatively small-scale systems \cite{haller18} and applications to large scale problems are still the object of intensive research \cite{NNM21}. Thus, the projection of dynamical systems onto low-dimensional nonlinear manifolds is still an open challenge. Data-driven approaches based on Proper Orthogonal Decomposition (POD) have been recently successfully applied to MEMS by Gobat et al.\cite{POD21}; however, these are still linear approaches that might experience some difficulties in reproducing the correct curvature of invariant manifolds while retaining only few DOFs in the ROM, and eventually fail to reach real-time performances.

Exploiting machine learning methods for constructing surrogates
of dynamical systems has recently been an
area of increasing interest for the system dynamics community.
Among others, huge success has been encountered by the Physics Informed
Neural Networks (PINN) \cite{Raissi19} which have been applied in multiple contexts and also in solid mechanics \cite{Raissi21}. These neural networks leverage the automatic differentiation feature to directly enforce partial differential equations. However, so far they still do not serve the aim of generating ROMs for parameter dependent problems.

Deep Learning (DL) techniques come as an inspiration to handle the complex reduction process of dynamical systems, unvealing low-dimensional features from black-box data streams \cite{Guo2018,Guo2019,Carlberg20}. 
Brunton and coworkers recently applied the SINDy method \cite{kaiser2018sparse,brunton2016discovering, brunton2016sparse} in combination with autoencoder neural networks to discover the underlying model of dynamical systems \cite{champion2019data}.
Fresca et al.\cite{FrescaManzoni20,FrescaManzoni21}, 
proposed a non-intrusive DL-based ROM technique, which we refer to as DL-ROM. 
Combining in a suitable way a convolutional autoencoder (AE) and a deep feedforward neural network (DFNN), the DL-ROM technique enables the construction of an efficient ROM, whose dimension is as close as possible to the number of parameters upon which the solution of the differential problem depends. The encoder part performs an operation of feature extraction forcing the high dimensional data to be reduced, at the bottleneck layer, to few latent variables. We highlight that this approach builds on the idea that the system dynamics develops on a very low dimensional curved invariant manifold. 
Indeed this is very often the case for real MEMS devices which, in the linear regime, oscillate according to a specific mode of interest, while this turns into the corresponding Nonlinear Normal Mode (NNM) as the actuation intensity increases \cite{touze2006,NNM21}. 
The procedure exploits snapshots taken from a set of FOM solutions (for selected parameter values and time instances) and deep neural network architectures to learn, in a non-intrusive way, (i) the nonlinear trial manifold where the ROM solution is sought, (ii) the nonlinear reduced dynamics, and (iii) and the reconstruction of the approximate FOM response starting from the latent variables. In the framework proposed by Fresca et al.\cite{FrescaManzoni21}, a further dimensionality reduction carried out on the FOM data through proper orthogonal decomposition (POD), yielding the so-called POD DL-ROM technique, also allows to speedup training times and to compress data dimensions, enhancing the construction of DL-ROMs. 
Neural networks have also been used to model simple structures by Simpson et al. \cite{Chatzi21}, with two main differences: 
the convolutional autoencoder uses Long Short-Term Memory (LSTM) networks; a statistical regression model is used instead of the DFNN. The adoption of LSTM cells is beneficial for the prediction of the time evolution beyond the training window, which is not targeted in the present contribution focusing on periodic responses. More in general, DL algorithms in conjunction with POD have already been exploited to address long-term predictions in time, however without including parameter dependencies \cite{wang2018model,lui_wolf_2019}. Recurrent neural networks have been considered by Gonzalez and Balajewicz\cite{Gonzalez18} to evolve low-dimensional states of unsteady flows, exploiting either POD or a convolutional recurrent autoencoder to extract low-dimensional features from snapshots. DL algorithms have also been used to describe the reduced trial manifold where the approximation is sought, then exploiting a minimum residual formulation to derive the ROM as done by Lee and Carlberg\cite{Carlberg20}.

In this paper we propose a new version of the DL-ROM discussed above, in which a first 
dimensionality reduction of the data is achieved 
by means of a POD-Galerkin (POD-G) ROM which has been recently applied to microstructures by Gobat et al.\cite{POD21} showing very good predictive capabilities.
In this way, the costly FOM data generation phase is greatly reduced and the training of the 
DL-ROM is performed using cheaper POD-G snapshots covering the whole parameter range. 
This allows addressing real industrial examples in complex scenarios where multiplicity of solutions exist thus requiring specific provisions.
Compared to other surrogate models exploiting machine/deep learning algorithms, a distinguishing feature of POD-G DL-ROMs is their capability to compute the whole solution field, for any new parameter instance and time instant, at testing time, thus enabling the extremely efficient evaluation of any output quantity of interest depending on the solution field.
The technique proposed in this work can be considered as the Data Driven counterpart
of the Direct Normal Form method for invariant manifolds that has been 
recently applied to large scale finite element systems of mechanical structures \cite{vizza2020,NNM21,vizzaccaro2021high}.

The structure of the paper is as follows.
After a brief definition in Section~\ref{sec:formulation} of the problem at hand and of the
type of nonlinearities involved,
we discuss the proposed data-driven procedure in Section~\ref{sec:ROM}.
The POD-G ROM, i.e.\ the first level of dimensionality reduction, is presented in Section~\ref{sec:PODROM}, while the second Deep Learning level is analysed in
Section~\ref{sec:PODGDLROM}, detailing the structure of the Autoencoder and of the FF Neural Network.
Two applications are then presented in Section~\ref{sec:applications},
starting with the academic example of a clamped-clamped beam in Section~\ref{sec:ccbeam},
and later addressing a real micromirror in Section~\ref{sec:micromirror}.

\section{Problem formulation}
\label{sec:formulation}

We will consider MEMS devices with given fixed geometry subjected to a periodic forcing of angular frequency $\omega$ and modulated in amplitude by the coefficient $\beta$. The input parameter vector is $\bfmu=[\omega,\beta] \in\mathcal{P}\subset \Rbb^{n_{\bfmu}}$, with $\mathcal{P}$ a closed and bounded set of dimension $n_{\bfmu}=2$.  
These devices are usually characterised by several Frequency Response Functions (FRFs). For instance, in the specific case of the doubly clamped beam analysed in Section~\ref{sec:ccbeam}, the FRF plotted in Figure~\ref{fig:beam_POM}f) expresses the mid-span displacement $u_{\tmid}$
of the beam in terms of $\bfmu$. Each dot represents the maximum value of
$u_{\tmid}$ in the corresponding periodic solution. 
Moreover one might desire to monitor the evolution of stresses or strains 
at other specific points in order to guarantee a correct operation of the device.
Providing all these output data is the main goal of any numerical technique.
Therefore we aim to formulate a fast and reliable {\it reduced-order} technique capable of predicting periodic responses of period $T$, 
with $T=2\pi/\omega$, and of providing as output the corresponding full field
response in terms of displacements, strains and stresses.


Many microsystems, like e.g.\ micromirrors, undergo large transformations and in particular large rotations.
The device initially occupies the domain $\Omega_0$ described by material coordinates $\bfX$
and is subjected to the transformation $\bfx=\bfphi(\bfX,t)=\bfX+\bfu(\bfX,t)$, with
$\bfu$ displacement vector and $\bfx$ actual spatial coordinates.
The boundary $\partial \Omega_0$ 
is partitioned in $\partial\Omega_D$ and  $\partial\Omega_N$ where Dirichlet and Neumann periodic boundary conditions are enforced, respectively. 
In this contribution, for the sake of simplicity, only homogeneous 
boundary conditions will be considered.
By assumption, this partition does not depend on the transformation thus ruling out specific applications including, e.g.\ evolving contact.
The device is actuated by given time-periodic body forces $\bfB(\bfX,t)$.
The overall governing system of PDEs with periodic time boundary conditions reads:
%
\begin{subequations} 
	\label{eq:strongmech}
	\begin{align}
		\rho_0 \ddot{\bfu}(\bfX,t) - \nabla\cdot\bfP(\bfX,t)- 
		\rho_0 \bfB(\bfX,t;\bfmu)=\bfzero \; 
		& \quad \text{for} \; (\bfX,t)\;\text{in} \; \Omega_0\times[0,T], \label{eq:strongmech_a}
		\\
		\bfP(\bfX,t)\cdot\bfN(\bfX)= \bfzero \; 
		& \quad \text{for} \; (\bfX,t)\;\text{in} \; \partial\Omega_N\times[0,T], \label{eq:strongmech_c}
		\\
		\bfu(\bfX,t)=\bfzero \; 
		& \quad \text{for} \; (\bfX,t)\;\text{in} \; \partial\Omega_D\times[0,T], \label{eq:strongmech_d}
		\\
		\bfu(\bfX,0)=\bfu(\bfX,T) \; 
		& \quad \text{for} \; \bfX\;\text{in} \; \Omega_0, \label{eq:strongmech_e}
		\\
		\dot{\bfu}(\bfX,0)=\dot{\bfu}(\bfX,T) \; 
		& \quad \text{for} \; \bfX\;\text{in} \; \Omega_0, \label{eq:strongmech_ee}
		\\
		\bfS(\bfX,t)=\cA(\bfX):\bfE(\bfX,t) \; 
		& \quad \text{for} \; (\bfX,t)\;\text{in} \; \Omega_0\times[0,T], \label{eq:strongmech_b}
        \\
        \bfE(\bfX,t) =\frac{1}{2}\left(\nabla\bfu(\bfX,t)+\nabla^T\bfu(\bfX,t)+\nabla^T\bfu(\bfX,t)
        \cdot\nabla\bfu(\bfX,t)\right)  
        &   \quad \text{for} \;  (\bfX,t)\;\text{in} \; \Omega_0\times[0,T], \label{eq:strongmech_gg}
        \\
        \bfP(\bfX,t) =(\tm+\nabla \bfu(\bfX,t))\cdot\bfS(\bfX,t),  \label{eq:strongmech_gh}
        &   \quad \text{for} \;  (\bfX,t)\;\text{in} \; \Omega_0\times[0,T].
	\end{align}
\end{subequations}
Eq.~\eqref{eq:strongmech_a} expresses the conservation of momentum
where $\rho_0$ is the initial density  and $\bfP$ is the first Piola-Kirchhoff stress
\cite{malvern}.
Eq.~\eqref{eq:strongmech_c} and \eqref{eq:strongmech_d} are the 
Neumann and Dirichlet boundary conditions respectively. 
Eqs.~\eqref{eq:strongmech_e} and 
\eqref{eq:strongmech_ee} enforce the periodicity condition on displacements
and velocities. 
The device is made of cubic single crystal silicon or polysilicon, thus admitting only
small strains,  a condition which is well described by the Saint Venant-Kirchhoff constitutive model Eq.~\eqref{eq:strongmech_b} 
between the Second Piola-Kirchhoff stress $\bfS$ and the Green-Lagrange strain tensor $\bfE$ 
(Eq.~\eqref{eq:strongmech_gg})
through the fourth-order elasticity tensor $\cA$ endowed with major and minor symmetries.
Finally, Eq.~\eqref{eq:strongmech_gh} formulates the link between the two Piola-Kirchhoff stress
tensors \cite{malvern}.


The weak form of the momentum conservation, i.e. Eqs.~\eqref{eq:strongmech_a}-\eqref{eq:strongmech_c}, reads:
\begin{align}
	\label{eq:PPV}
		\int_{\Omega_0}\rho_0\ddot{\bfu}(\bfX,t)\shcdot\bfw(\bfX)\,\dOmega_0  +\!
		\int_{\Omega_0}\bfP(\bfX,t)\shddot\nabla^T\!\bfw(\bfX)\,\dOmega_0  = 
		 \int_{\Omega_0}\!\rho_0 \bfB(\bfX,t;\bfmu) \shcdot\bfw(\bfX) \,\dOmega_0
	 \qquad \forall \bfw\in H^1_0(\Omega_0),
\end{align}
where $\bfw$ is the test velocity selected in $H^1_0(\Omega_0)$, i.e.\ the space of functions with finite energy that vanish on the portion $\partial\Omega_D \subset \partial\Omega_0$ where Dirichlet boundary conditions are prescribed.
Within the present context, it is worth stressing that Eq.~\eqref{eq:PPV} exactly accounts for geometric (elastic and inertia) nonlinearities, e.g., large rotations or nonlinear mode coupling.  

The spatial discretization of Eq.~\eqref{eq:PPV}, e.g.\ by means of the finite element method, with the additional inclusion of a Rayleigh model damping term, yields a system of coupled nonlinear differential equations representing the full order model (FOM):
\begin{subequations} 
	\label{eq:PPV_d2}
	\begin{align}
		&\bfM \ddot{\bfu}_h(t) + \bfC\dot{\bfu}_h(t)  + \bfK\bfu_h(t) + \bfG(\bfu_h,\bfu_h) + \bfH(\bfu_h,\bfu_h,\bfu_h) = 
		\bfF(t;\bfmu), \qquad t \in (0,T)\\
		&\bfu_h(0)=\bfu_h(T), \qquad
		\dot{\bfu}_h(0)=\dot{\bfu}_h(T) 
	\end{align}
\end{subequations}
where the vector $\bfu_h(t) \in \mathbb{R}^{N_h}$ collects the $N_h$ unknown displacements nodal values, 
$\bfM \in \mathbb{R}^{N_h \times N_h}$ is the mass matrix, $\bfC=(\omega_0/Q)\bfM$ is the Rayleigh mass-proportional 
damping matrix -- considering a reference eigenfrequency $\omega_0$ and a quality factor $Q$ -- and  
$\bfF(t;\bfmu) \in \mathbb{R}^{N_h}$ is the nodal force vector which depends on the vector of parameters $\bfmu$. 
The internal force vector has been exactly decomposed in linear, quadratic, and cubic power terms of the displacement: $\bfK \in \mathbb{R}^{N_h\times N_h} $ is the stiffness  matrix related to the linearized system, while $\bfG \in \mathbb{R}^{N_h}$ and $\bfH \in \mathbb{R}^{N_h}$ are  vectors given by monomials of second and third order, respectively. 
We stress that the components of these vectors can be expressed using an indicial notation as
\[
G_i=\sum_{j,k=1}^{N_h}G_{ijk} u_{h,j} u_{h,k}, \quad   H_i=\sum_{j,k,l=1}^{N_h}H_{ijkl}u_{h,j} u_{h,k} u_{h,l}, \qquad i=1,\ldots, N_h.
\] 
Eq.~\eqref{eq:PPV_d2} represents our high-fidelity FOM which depends on the input parameters
$\bfmu$. 
Our goal is the efficient numerical approximation of the solution manifold:
\begin{equation}
	\label{eq:sethifi}
	\cS=\left\{\bfu_h(t;\bfmu) |\; t \in [0,T)\; \text{and}\; \bfmu \in \cP \subset \Rbb^{n_{\bfmu}} \right\}\subset \Rbb^{N_h}.
\end{equation}
In particular, assuming that, for any given parameter $\bfmu \in \mathcal{P}$, the FOM admits a unique solution  for each $t \in (0,T)$, that is,  the dynamical system Eq.~\eqref{eq:PPV_d2} has a unique trajectory for each parameter instance, the intrinsic dimension of the solution manifold is at most $n_{\bfmu} + 1  \ll N_h$, where $n_{\boldsymbol \mu}$ is the number of parameters (time thus plays the role of an additional coordinate).

The numerical solution of the FOM Eq.~\eqref{eq:PPV_d2} to compute the steady state response 
is a challenge in itself for large scale problems.
One option is the use of time marching methods to simulate a sufficiently large number $N_c$ of cycles, 
where $N_c$ is typically inversely proportional to the damping. 
This technique resorts to robust algorithms implemented in most of the commercial software, but when damping is very small, as in most MEMS devices, the computational effort may not be affordable. 
In other approaches, like the Harmonic Balance (HB) method \cite{actuators21,KerschenHB2015}, the 
unknown displacements are expressed as the sum of Fourier components thus automatically respecting periodicity conditions. However, their implementation requires dedicated codes and non-standard computing facilities.

\section{Data-driven Model Order Reduction}
\label{sec:ROM}

In this section, we briefly review the construction of POD-Galerkin (POD-G)  ROMs and of the POD-G Deep Learning-based ROM (POD-G DL-ROM) technique which combines the best features of POD and DL algorithms in order to achieve real-time simulations of the structural behavior of MEMS devices.

\subsection{POD-Galerkin ROM (POD-G ROM)}
\label{sec:PODROM}

A simple way to build a ROM able to approximate Eq.~\eqref{eq:sethifi} is through the introduction of a linear reduced trial manifold $\tilde{\cS}_h^{N, lin}=\textnormal{Col}(\mathbf{V})$ of dimension $N\lll N_h$, spanned by the $N$ columns of a matrix $\mathbf{V} \in \Rbb^{N_h\times N}$. Thus we can achieve an approximation $\tilde{\bfu}_h(t,\bfmu)\approx\bfu_h(t,\bfmu)$ as:
\begin{equation}
	\tilde{\bfu}_h(t,\bfmu)=\mathbf{V}\bfu_N(t;\bfmu)
	\label{eq:trialspace}
\end{equation}
in which  $\tilde{\bfu}_h:[0,T)\times\cP\rightarrow\tilde{\cS}_n$ and $\bfu_N(t;\bfmu)\in\Rbb^N$, $\forall t\in [0,T)$ and $\bfmu\in \cP$, is the vector of degrees of freedom of the ROM approximation.
The most commonly used technique to build the subspace $\tilde{\cS}_h^N$ bases is through POD.  
The first step in the construction of a POD-G ROM requires to generate 
a  matrix $\bfS_u\in \mathbb{R}^{N_h\times N_s}$, whose $N_s$ columns collect snapshots 
of the FOM solutions, obtained for different values of the parameters $\bfmu$
\begin{equation}
	\bfS_u=[\bfu_N(t^1,;\bfmu_1) | ...|\bfu_N(t^{N_t},;\bfmu_1)|...|\bfu_N(t^1,;\bfmu_{N_{\ttrain}}) | ...|\bfu_N(t^{N_t},;\bfmu_{N_{\ttrain}})].
\end{equation}
Each snapshot in $\bfS_u$ is sampled inside a certain time interval $[0,T]$ partitioned in $N_t$ time steps, with $t^k, k=1, \ldots, N_t$, and the parameter space $\mathcal{P}$ where a set of $N_{\ttrain}$ parameter instances $\bfmu_k,k=1,\ldots,N_{\ttrain}$, is selected.
Next, POD computes the  Singular Value Decomposition (SVD) of the matrix $\bfS_u$ 
\[
\bfS_u = \bfU \bfSigma \bfZ^T,
\]
where the columns of the orthonormal matrix $\bfU\in \mathbb{R}^{N_h\times N_h}$ are 
the left singular vectors and
the columns of the orthonormal matrix $\bfZ\in \mathbb{R}^{N_s\times N_s}$ are the right singular vectors.
The diagonal elements of $\bfSigma\in \mathbb{R}^{N_h\times N_s}$ are the singular values of the matrix $\bfS_u$ and are conventionally ordered
from the largest to the smallest. In particular, POD selects the columns of $\mathbf{V}$ in Eq.\eqref{eq:trialspace} as the first $N$ left singular vectors of $\bfS_u$, often called POD modes (POMs) in the literature 
\cite{kerschen2005,alfiobook,lu2019}.

Once the linear trial POD subspace has been obtained, projecting the FOM Eq.~\eqref{eq:PPV_d2} onto the POD subspace yields the structural dynamics geometric POD-G ROM, under the form of a $N$-dimensional nonlinear ODE system, whose solution provides the dynamics of the generalized coordinates
$\bfu_N\in \mathbb{R}^{N}$: 
\begin{equation}\label{eq:PPV_POD}
		\bfM\aPOD \ddot{\bfu}_N + \bfC\aPOD\dot{\bfu}_N  + \bfK\aPOD\bfu_N + \bfG\aPOD(\bfu_N,\bfu_N) 
		 +\bfH\aPOD(\bfu_N,\bfu_N,\bfu_N) = \bfF\aPOD(\bfu_N,t;\bfmu), \quad t \in (0,T),
\end{equation}
where
\[
\bfM\aPOD=\bfV^T\bfM\bfV, \quad \bfC\aPOD=\bfV^T\bfC\bfV, \quad \bfK\aPOD=\bfV^T\bfK\bfV,
\]
\[
\bfF\aPOD=\bfV^T\bfF, \quad  G\aPOD_i=g\aPOD_{ijk} u_{N,j} u_{N,k}, \quad  H\aPOD_i=h\aPOD_{ijkl} u_{N,j} u_{N,k} u_{N,l},
\] 
with $\bfM\aPOD,\bfC\aPOD,\bfK\aPOD \in \mathbb{R}^{N\!\times\!N}$.
The computation of the vectors $\bfG\aPOD$ and $\bfH\aPOD$ entails $O(N^3)$ and $O(N^4)$ terms, respectively.
Note that the coefficients $g\aPOD_{ijk}$ and $h\aPOD_{ijkl}$ can be precomputed, and that the reduced problem can be assembled efficiently thanks to its polynomial nature, 
thus avoiding the use of  hyper-reduction techniques such as the (discrete) 
empirical interpolation 
method \cite{barrault2004anempirical,chaturantabut2010nonlinear,maday2008ageneral}.

The POD-G approach has been recently benchmarked on several MEMS including beams, arches and mirrors. With reference to MEMS structures Gobat et al.\cite{POD21}  have shown that, 
provided that a sufficient number of POD modes is included in the trial space, the technique accurately reproduces the response of the device. However, in particular when large rotations are involved and nonlinearities are strong, the dimension of the trial space increases and solution of the ROM with continuation techniques comes with a high computational cost,
failing to serve the final goal of generating a real-time simulation tool.

\subsection{POD-Galerkin-enhanced deep learning based-reduced order models (POD-G DL-ROMs)}
\label{sec:PODGDLROM}

POD-DL-ROMs are non-intrusive ROMs, which aim at approximating the map $(t, \bfmu) \rightarrow {\bf u}_h(t, \bfmu)$ by describing both the trial manifold and the reduced dynamics through deep neural networks \cite{fresca2020POD}. To reduce the dimensionality of the snapshots and avoid feeding training data of very large dimension $N_h$, POD is first applied -- realized through randomized SVD (rSVD) \cite{halko2011finding} -- to the snapshot set $\bfS_u$; then, a DL-ROM is built to approximate the map between $(t, \bfmu)$ and the POD generalized coordinates $\bfu_N(t; \bfmu)$. 

The POD-DL-ROM technique is extremely efficient and it is able to model highly nonlinear problems by identifying the manifold underlying the dynamics of the system in a complete data-driven and black-box, non-intrusive way. On the other hand, its data-driven nature implies that the sampled data provided must span the parameter space of interest and contain all the information necessary to accurately approximate the solution manifold. In this way, the size of the training dataset increases with the number of parameters considered in the FOM. 

Instead, POD projection-based ROMs exploit a data-driven process to generate the basis of the linear subspace, later this basis is used to project the nonlinear system given by the spatial discretisation. 
For instance (exact) projection procedures like the POD-G ROM outlined in Section~\ref{sec:PODROM}, lead to a nonlinear system of ODEs whose degrees of freedom can be physically interpreted and that can be solved with many different integration packages (i.e. Auto07p \cite{Doedel}, MANLAB\cite{Guillot1}, Nlvib \cite{NVLIB}, Matlab built-in functions, BifurcationKit \cite{BifurcationKit} etc.). On the other hand, such POD-G ROM is strongly intrusive and it is way less efficient than a DL-ROM, because the ODE system must be solved up to a given time instant of interest for a given parameter instance. These conditions may be acceptable when a limited number of inquiries are needed, but when a refined inspection of the system response dynamics is necessary for optimal design purpose or when a real-time solution is needed (e.g. online control) POD-G ROMs, and in general intrusive approaches, fail. 

\begin{figure}[htb]
	\centering
	\includegraphics[width =.85\textwidth]{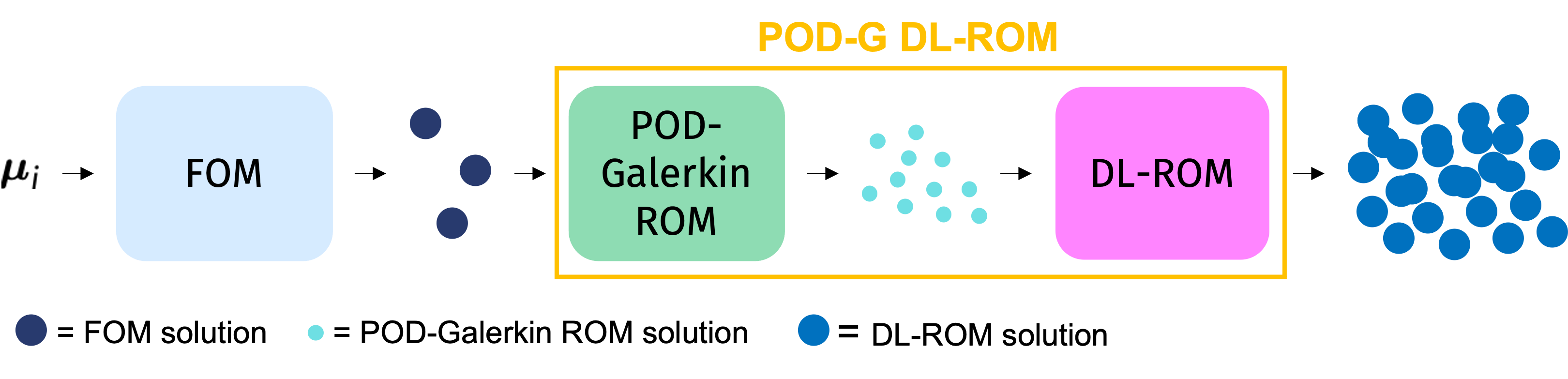}
	\caption{POD-G DL-ROM procedure scheme. By solving the FOM for few parameter instances, a small number of FOM snapshots are generated. The high-fidelity snapshots are used to compute the POD basis matrix $\bfV$ and generate the POD-G ROM. This is later exploited to compute a large number of reduced snapshots, for different parameter instances, able to finely span the parameter space of interest. These snapshots represent the input to the POD-G DL-ROM neural network. Once the POD-G DL-ROM is trained, it is used to generate all the approximate solutions required.}
	\label{fig:POD-G DL-ROM_scheme}
\end{figure}

In this contribution, we exploit the most appealing features of both techniques to achieve a real-time simulation of a MEMS. 
A sketch of the proposed  procedure, named POD-Galerkin-enhanced Deep Learning-based ROM (POD-G DL-ROM), is outlined in Fig.~\ref{fig:POD-G DL-ROM_scheme}.

First, a few FOM simulations, that loosely span the parameter space, are performed. The generated high-fidelity snapshots are processed with SVD and a subset of POD modes is selected. The number of POD modes kept in the subspace must be defined in order to guarantee a good approximation of the underlying manifold as shown by Gobat et al.\cite{POD21}, where POD-G ROMs are applied to mechanical systems with low damping like MEMS. By exploiting the generated linear subspace, the POD-G ROM is built following the procedure detailed in Section~\ref{sec:PODROM}. The POD-G ROM depends explicitly on the model parameters of interest, i.e. the load amplitude $\beta$ and the forcing frequency $\omega$, and it is solved, at testing time, to create snapshots of the ROM intrinsic coordinates finely spanning the parameter space $\mathcal{P}$. 

The resulting reduced solutions are provided as input to the POD-G DL-ROM neural network during the training stage. 
Following the DL-ROM technique outlined by Fresca et al.\cite{FrescaManzoni21}, the DFNN is simultaneously trained to provide as output the minimal coordinate vector for the same parameter instance. Finally, the decoder function of the convolutional AE maps the minimal coordinates to the approximated intrinsic coordinates. 

More precisely, the POD-G DL-ROM approximation of the FOM solution ${\bfu}_h(t; \bfmu)$ is  
\begin{equation*}
\tilde{\bfu}_h(t;\bfmu,\bftheta_{DF},\bftheta_{D}) = \mathbf{V} \tilde{\bfu}_N(t; \bfmu, \bftheta_{DF}, \bftheta_{D}) \approx \mathbf{V} {\bfu}_N(t; \bfmu, \bftheta_{DF}, \bftheta_{D}) \approx \bfu_h(t; \bfmu),
\end{equation*}
that is, it is sought in a linear trial manifold of (potentially large) dimension $N$, 
\begin{equation}
\label{manifold_linear_N}
\tilde{\mathcal{S}}_h^{N} = \{
{\bf V} \tilde{\bfu}_N(t; \bfmu, \bftheta_{DF}, \bftheta_{D})  \ |  \; \tilde{\bfu}_N(t; \bfmu, \bftheta_{DF}, \bftheta_{D}) \in \mathbb{R}^N,    \ t \in [0, T)  \; , \bfmu \in \mathcal{P}   \} \subset \mathbb{R}^{N_h},
\end{equation}
by applying the DL-ROM strategy \cite{FrescaManzoni21} to approximate ${\bfu}_N(t; \bfmu)$  -- rather than directly $\mathbf{V}^T\bfu_{h}(t ; \bfmu )$.  The DL-ROM approximation 
$\tilde{\bfu}_N(t; \bfmu, \bftheta_{DF}, \bftheta_{D}) \approx {\bfu}_N(t ; \bfmu )$ takes the form
\begin{equation}
\tilde{\bfu}_N(t; \boldsymbol \mu, \bftheta_{DF}, \bftheta_D ) = {\mathbf{f}}^D_N({\bfphi}_n^{DF}(t; \bfmu, {{\bftheta}_{DF}}); \bftheta_{D}),
\label{eq:u_N_approx}
\end{equation}
and is sought in a reduced nonlinear trial manifold $\tilde{\mathcal{S}}_N^n$ of very small dimension $n \ll N$; usually, $n \approx n_{\bfmu} + 1$ -- here time is considered as an additional parameter. As in a POD-DL-ROM, in a POD-G DL-ROM both the reduced dynamics and the reduced nonlinear manifold (or trial manifold) where the ROM solution is sought must be learnt. In particular:
\begin{itemize}
\item \textit{Reduced dynamics learning.} To describe the system dynamics on  the nonlinear trial manifold $\tilde{\mathcal{S}}_N^n$,  the intrinsic coordinates of the approximation $\tilde{\bfu}_N$ are defined as  
\begin{equation*}
\label{eq:phi_n}
\bfu_n(t; \boldsymbol \mu) = \bfphi_n^{DF}(t; \boldsymbol \mu, \bftheta_{DF}),
\end{equation*}
where   $\bfphi^{DF}_n(\cdot ;  \cdot, \bftheta_{DF}) : [0, T) \times \mathbb{R}^{n_{\mu}} \rightarrow \mathbb{R}^n$ is a DFNN, 
consisting of the repeated composition of a nonlinear activation function,  applied to a linear transformation of the input, multiple times.  Here $\bftheta_{DF}$ denotes the  DFNN {parameters} vector, collecting the weights and biases of each of its layers;
\item \textit{Nonlinear trial manifold learning.} To model the  reduced nonlinear trial manifold $\tilde{\mathcal{S}}_N^n$, we  employ the decoder function of a convolutional autoencoder (AE), that is, 
\begin{equation}
\tilde{\mathcal{S}}_N^n = \{ \tilde{\bfu}_{N}(t; \boldsymbol \mu) =  {\mathbf{f}}^D_N(\bfphi_n^{DF}(t; \bfmu, {{\bftheta}_{DF}}); \bftheta_{D}) \; | 
 \; \bfphi_n^{DF}(t; \bfmu,  {{\bftheta}_{DF}}) \in \mathbb{R}^{n}, \ t \in [0, T) \; , \; \bfmu \in \mathcal{P} \subset \mathbb{R}^{n_{\mu}} \} 	\subset \mathbb{R}^N,
\label{eq:manifold_tilde_S_n}
\end{equation}
where ${\mathbf{f}}^D_N( \cdot ; \bftheta_{D}) : 	\mathbb{R}^n \rightarrow \mathbb{R}^N$ denotes the   decoder function of a convolutional AE obtained as   the composition of several layers (some of which are convolutional), depending upon a vector ${\bftheta}_{D}$ collecting all the corresponding weights and biases.
\end{itemize}
Finally, the encoder function ${\mathbf{f}}^E_n( \cdot ; \bftheta_{E}) : 	\mathbb{R}^N \rightarrow \mathbb{R}^n$  --  depending upon a vector ${\bftheta}_E$ of parameters -- of the convolutional AE can be used to map the POD-G reduced solution (or approximated intrinsic coordinates) $\bfu_{N} (t, \bfmu)$ associated to $(t, \bfmu)$ onto a low-dimensional representation
\begin{equation*}
\label{eq:f_n}
 {\tilde{\bfu}_n}(t; \bfmu, \bftheta_{E}) = {\mathbf{f}}_{n}^E({\bfu}_N(t; \bfmu); \bftheta_{E}).
\end{equation*}

\begin{figure}[htb]
	\centering
	\includegraphics[width =.85\textwidth]{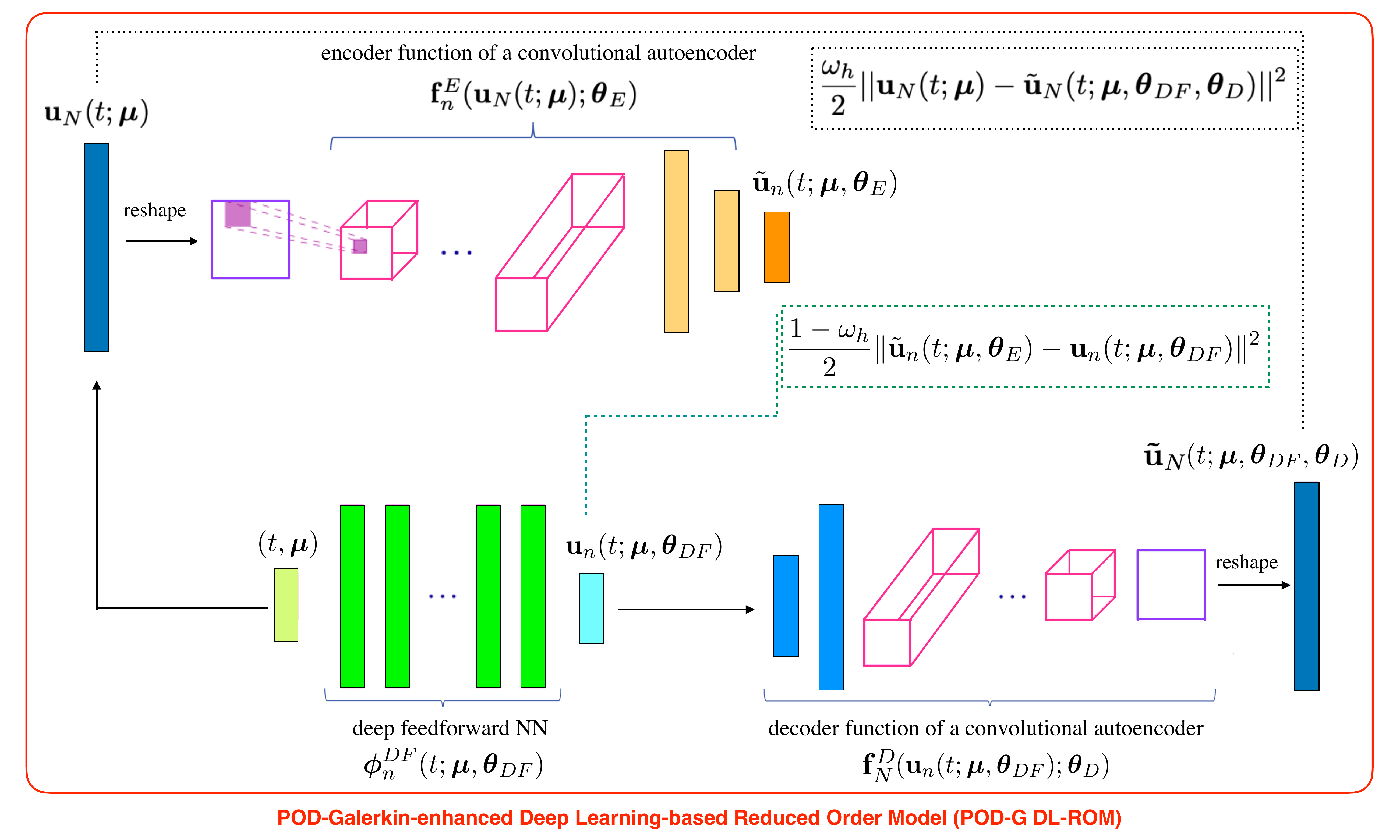}
	\caption{POD-G DL-ROM architecture. The approximated intrinsic coordinates $ \mathbf{u}_N(t; \boldsymbol{\mu})$ are provided as input to the encoder function which outputs $\tilde{\mathbf{u}}_n(t; \boldsymbol{\mu})$. The same parameter instance associated to the POD-G ROM reduced solution, i.e. $(t; \boldsymbol{\mu})$, enters the DFNN which provides as output $\mathbf{u}_n(t; \boldsymbol{\mu})$, and the error between the low-dimensional vectors is accumulated. The minimal coordinates $\mathbf{u}_n(t; \boldsymbol{\mu})$ are given as input to the decoder function returning the  approximated reduced solution $\mathbf{\tilde{u}}_N(t; \boldsymbol \mu)$. Then, the reconstruction error is computed. The approximation of the FOM solution is then recovered as $\tilde{\bfu}_h=\bfV\tilde{\bfu}_N(t;\bfmu)$.}
	\label{fig:POD-G DL-ROM}
\end{figure}

Hence, provided the parameter matrix $\mathbf{L} \in \mathbb{R}^{(n_{\boldsymbol \mu} + 1) \times N_{s}}$ defined as
\begin{equation}
\mathbf{L} = [(t^1, \bfmu_1) | \ldots | (t^{N_t}, \bfmu_1) | \ldots | (t^1, \bfmu_{N_{\ttrain}}) | \ldots | (t^{N_t}, \bfmu_{N_{\ttrain}})],
\end{equation}
and the snapshot matrix $\bfS_u$, training a POD-G DL-ROM requires to solve the optimization problem 	 
 \begin{equation}
\min_{\bftheta} \mathcal{J}(\bftheta) = \min_{\bftheta} \frac{1}{N_{\ttrain} N_t}\sum_{i=1}^{N_{\ttrain}} \sum_{k=1}^{N_t} \mathcal{L}(t^k, \boldsymbol \mu_i; \bftheta),  \label{eq:minimization_problem}
\end{equation}
where the \textit{per-example} loss function $\mathcal{L}(t^k, \bfmu_i;  {\bftheta})$ is given by the sum of two terms
\begin{equation}
\label{eq:loss_N}
\begin{split}
\mathcal{L}(t^k, \bfmu_i;  {\bftheta})  =  \frac{\omega_h}{2} \| {\bfu}_N(t^k; \bfmu_i) & - \tilde{\bfu}_N(t^k; \bfmu_i,  {\bftheta_{DF}, \bftheta_D})\|^2  \\ & + \frac{1-\omega_h}{2}   
\| \tilde{\bfu}_n(t^k; \bfmu_i, \bftheta_E) -  {\bfu}_n(t^k; \bfmu_i,  {\bftheta_{DF}})\|^2,
\end{split}
\end{equation}
with $\bftheta = (\bftheta_E, \bftheta_{DF}, \bftheta_D)$. The former term in Eq.~\eqref{eq:loss_N} is the reconstruction error  between the POD-G ROM reduced solutions and the POD-G DL-ROM approximations while the latter is the misfit between the  {intrinsic coordinates} and the output of the encoder function. Finally, $\omega_h \in [0,1]$ is a prescribed weighting parameter. The architecture of the POD-G DL-ROM neural network is shown in Figure~\ref{fig:POD-G DL-ROM}; note that, at testing time, we can discard the encoder function of the convolutional AE.


We remark that the computation of the low-fidelity snapshots, obtained as reduced solution of a POD-G ROM for several parameter instances, can be considered as an extension of the idea at the basis of the POD-DL-ROM technique, where the high-fidelity, FOM snapshots are projected onto the reduced linear trial subspace, generated by means of rSVD, before being fed as input to the AE. The main advantage is that solving the POD-G ROM is more efficient than the solution of the FOM, thus we can use it to generate all the snapshots required by the POD-G DL-ROM neural network to be able to generalize to new, unseen scenarios by remarkably decreasing the time required for the generation of the training/testing datasets. On the other hand, the quality of the reduced solution is influenced by the number of retained POD modes, thus the preparation of the POD-G ROM is crucial.

Summarizing the procedure, the POD-G DL-ROM exploits the POD-G ROM in the training stage to: 1) finely sample the parameter space to guarantee the proper generalization capabilities of the DL-ROM neural network; 2) perform a dimensionality reduction of the FOM with consequent advantages in memory management and efficiency during the training stage of the DL-ROM; 3) compute the POD basis matrix employed to recover $\bfu_h$. 
The DL-ROM is still used to approximate the FOM solution during the online, testing stage, thus its efficiency and generalization capabilities are fully exploited. As result, we get a real-time simulation of the device dynamics, whose performance is detailed in 
Section~\ref{sec:applications}.


\subsection{Frequency vs phase control}
\label{sec:omegaphi}

As recalled in Section~\ref{sec:intro}, we aim at simulating the FRF of MEMS.
However, the intrinsic nonlinearities imply that the system response may not be uniquely defined for a given external frequency $\omega$ and load $\beta$ as a consequence of bifurcations, e.g. saddle-node bifurcations. The steady-state response of the system will converge to a certain branch depending on the initial conditions.  
This is the case, for instance, of the simple clamped-clamped beam analysed in Section~\ref{sec:ccbeam}
which is representative of a large class of resonating MEMS and it is characterized by a FRF
like the one in Figure~\ref{fig:beam_POM}.
For a given load multiplier $\beta$ and angular frequency $\omega$,
the maximum midspan displacement over the period $T = 1/\omega$ of the beam can have one or three amplitude values, due to the so-called
hardening effect.
This portrait of the dynamics may be further complicated in presence of internal resonances where the number of overlapping branches and alternative solutions increases.


\begin{figure}[ht]
	\centering
	\includegraphics[width = .85\linewidth]{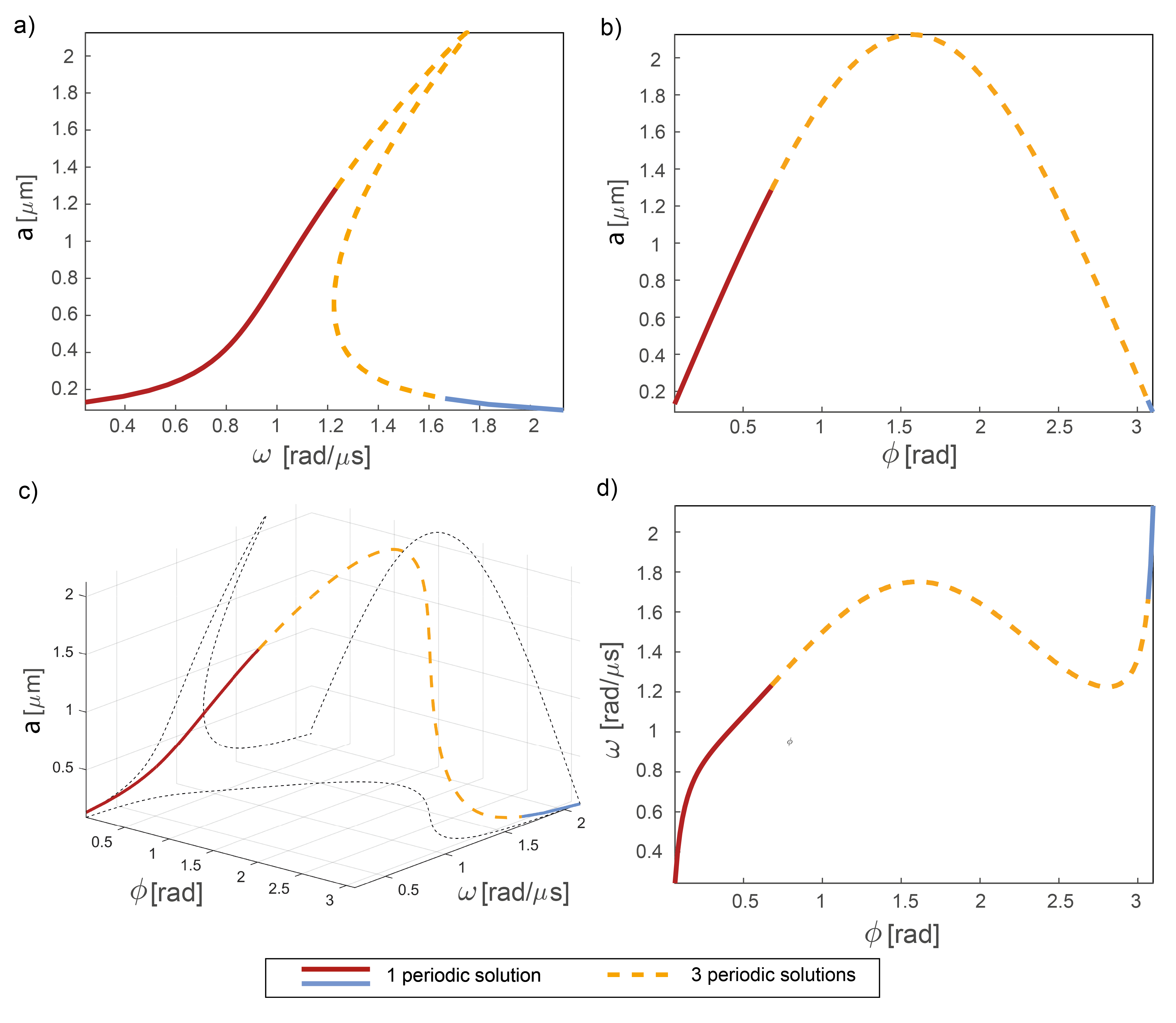}
	\caption{Response of the Duffing oscillator eq.\eqref{eq:duff} 
	for $\omega_0$=1[rad/$\mu$s], $\varepsilon$=1, $\beta=0.2$, $k_3$=0.5 and $\xi$=0.05. Figure a): typical FRF, where $a$ is plotted versus the forcing angular frequency $\omega$. Three regions can be identified: two 
	regions (red and blue lines, respectively) characterised by a single periodic solution for a given $\omega$, and a region (orange line) where three solutions exist. This plot can be viewed as the projection on the $\omega-a$ plane of the response in the 3D space $\phi-\omega-a$ (Figure c).
    Figure b): projection of the system response of Figure c) on the $\phi-a$ plane, where 
    uniqueness is preserved for any phase value. 
    Figure d): projection of the system response of Figure c) on the $\phi-\omega$ plane.}
	\label{fig:duff_phase}
\end{figure}

In this contribution, we will focus on systems without internal resonances and forced with a single frequency excitation close to an eigenfrequency. In these conditions, the dynamics is similar to the one of a simple Duffing resonator:
%
\begin{equation} 
\label{eq:duff}
	\ddot{u}+\omega_0^2 u+\varepsilon( 2\xi \omega_0 \dot{u} +k_3 u^3)=\varepsilon \beta \cos(\omega t),
\end{equation}
where $\omega_0$ is the eigenfrequency, $\xi$ is damping coefficient, $k_3$ cubic nonlinearity coefficient, $\beta$ is the load multiplier and $\varepsilon\ll 1$ is a {\it bookkeeping} parameter that sorts the order of magnitude of the terms \cite{thomsen2003vibrations}.
The Multiple Scales solution of Eq.~\eqref{eq:duff} reads:
\begin{equation} \label{eq:duff_sol}
	u(t)=a(\beta) \cos(\omega t-\phi)+\varepsilon a^3(\beta)\frac{k_3}{32\omega_0^2}\cos(3(\omega t-\phi)) + O(\varepsilon^3),
\end{equation}
%
%
%
%
which is plotted in several forms in Figure~\ref{fig:duff_phase}.
In particular, we remark that the amplitude $a(\beta)$ can be expressed as a single-valued function of the phase $\phi\in[0:\pi]$. 
It is worth recalling that this property of Duffing-like systems is 
exploited  experimentally in phase-controlled closed-loop experiments that allow tracking both stable and unstable branches \cite{jmems20reso,mestrom2009phase}. 

Since the POD modes in $\mathbf{V}$, eq.\eqref{eq:trialspace},
are energetically ordered and we know that the first POM is generally associated with the excited eigenmode, the corresponding dof in the POD-G ROM can be used to define the phase $\phi$.
Moreover, the periodic solution is expressed as the sum of several
harmonic components, each with a specific phase value, see e.g.\ Eq.~\eqref{eq:duff_sol}.
We hence define the phase referring to the lowest harmonic component which is typically the largest and also completely characterises linear problems.

Consequently, our problem will be no longer be expressed directly as a function of the parameters $\bfmu=[\beta,\omega]$, but rather of $\bfmu=[\beta,\phi]$.
However, since we aim at plotting FRFs, the relationship between $\phi$ and $\omega$ must be suitably approximated for each value of $\beta$. Coherently with the approach taken in this investigation, we use a DFNN to model the relationship  $(\phi, \beta) \mapsto \omega(\phi,\beta)$. The training data of the DFNN are represented by all the $\omega$ values at which the POD-G ROM solutions are sampled to feed the DL-ROM.
This way to proceed allows solving the lack of uniqueness of the response, but nevertheless the $\omega(\phi,\beta)$ function has two asymptotes for $\phi=0, \pi$, thus the regions close to these points are difficult to model with a DFNN. The best results have been obtained applying the following strategy.
Instead of using directly $\phi$ in the DFNN, we formulate the input of the
neural network in terms of $\sin\phi$ and $\cos\phi$. Then, the goal is to learn the function
\begin{equation}
f_{\phi-\omega}^{DF}  : \mathbb{R}^3 \rightarrow \mathbb{R} \quad \textnormal{such that} \quad (\sin\phi, \cos \phi, \beta) \mapsto f_{\phi-\omega}^{DF} \approx \omega(\phi,\beta).
\end{equation}
This simple feature engineering, based on the qualitative assumptions on the problem considered, allows mitigating the effects close to the boundaries, 
as shown in the following section.

\section{Numerical results}
\label{sec:applications}

We discuss the application of the proposed POD-G DL-ROM technique to MEMS devices: a clamped-clamped beam resonator and a micromirror. 
The FOM is solved through the HB method 
implemented, together with a continuation scheme, in a custom Fortran library \cite{actuators21} for large scale problems.
The generation of POD-G ROMs for MEMS has been extensively discussed by Gobat et al.\cite{POD21}. Here the ROM is solved with the Matlab package MANLAB \cite{Guillot1} that exploits the HB method and asymptotic expansion with continuation of periodic orbits.

\begin{table}[ht!]
\begin{center}
\begin{tabular}{|c|c|c|c|c|c|c|}
\hline
layer & input & output  & kernel  & $\#$of filters & stride & padding \\
  & dimension &   dimension &   size &  &  &   \\
\hline
1 & [N, N, 1] & [N, N, 8] & [5, 5] & 8 & 1 & SAME \\
\hline
2 & [N, N, 8] & [N/2, N/2, 16] & [5, 5] & 16 & 2 & SAME \\
\hline
3 & [N/2, N/2, 16] & [N/4, N4, 32] & [5, 5] & 32 & 2 & SAME \\
\hline
4 & [N/4, N/4, 32] & [N/8, N/8, 64] &[5, 5] & 64 & 2 & SAME \\
\hline
5 & N & 64 & & & & \\
\hline
6 & 64 & $n$ & & & &\\
\hline
\end{tabular}
\end{center}
\caption{Attributes of convolutional  and dense layers in the encoder $\mathbf{f}_n^E$.}
\label{table_transposed_convolutional_layers_encoder}
\end{table}

\begin{table}[ht!]
\begin{center}
\begin{tabular}{|c|c|c|c|c|c|c|}
\hline
layer & input & output  & kernel  & $\#$of filters & stride & padding \\
  & dimension &   dimension &   size &  &  &   \\
  \hline
1 & $n$ & 256 & & & &\\
\hline
2 & 256 & $N_h$ & & & &\\
\hline
3 & [N/8, N/8, 64] & [N/4, N/4, 64] & [5, 5] & 64 & 2 & SAME \\
\hline
4 & [N/4, N/4, 64] & [N/2, N/2, 32] & [5, 5] & 32 & 2 & SAME \\
\hline
5 & [N/2, N/2, 32] & [N, N, 16] & [5, 5] & 16 & 2 & SAME \\
\hline
6 & [N, N, 16] & [N, N, 1] & [5, 5] & 1 & 1 & SAME \\
\hline
\end{tabular}
\end{center}
\caption{Attributes of dense  and transposed convolutional layers in the decoder $\mathbf{f}_N^D$.}
\label{table_transposed_convolutional_layers}
\end{table} 

\subsection{Settings of neural networks}

Data normalization and  standardization enhance the training phase of the network by rescaling all the dataset values  to a common frame. For this reason, the inputs and the output of DL-ROM are  rescaled in the range $[0, 1]$ by applying an affine transformation. In particular, provided the parameter matrix $\mathbf{L} = [\mathbf{L}^{\ttrain}, \mathbf{L}^{\textnormal{val}}] \in \mathbb{R}^{(n_{\bfmu} + 1) \times N_s}$, we first rescale the time variable to the range $[0, 2\pi]$ as $\tau = t \omega$, where $\tau$ is the dimensionless time variable, and then we define
\begin{equation*}
\label{max_min_M}
L_{\tmax}^i = \max_{j =1, \ldots, N_s} L_{ij}^{\ttrain},  \qquad 
L_{\tmin}^i = \min_{j =1, \ldots, N_s} L_{ij}^{\ttrain}, \qquad i = 1, \ldots, n_{\bfmu} + 1,
\end{equation*}
so that   data are normalized by applying the following transformation
\begin{equation}
\label{eq:normalization}
L_{ij}^{\ttrain} \mapsto \frac{L_{ij}^{\ttrain} - L_{\tmin}^i}{L_{\tmax}^i - L_{\tmin}^i}, \qquad
i = 1, \ldots, n_{\bfmu} + 1,  \ j = 1, \ldots, N_s.
\end{equation} 
Each feature of the training parameter matrix is rescaled according to its maximum and minimum values. Regarding instead the training snapshot matrix $\mathbf{S}_u = [\mathbf{S}_u^{\ttrain}, \mathbf{S}_u^{\textnormal{val}}] \in \mathbb{R}^{N_h \times N_s}$, we define 
\begin{equation*}
\label{max_min_S}
S_{u,\tmax} = \max_{i =1, \ldots, N_h} \max_{j = 1, \ldots, N_s} S_{u,ij}^{\ttrain}, \quad \quad 
S_{u,\tmin} = \min_{i =1, \ldots, N_h} \min_{j = 1, \ldots, N_s} S_{u,ij}^{\ttrain}, 
\end{equation*}
and apply transformation Eq.~\eqref{eq:normalization}  by replacing $L_{\tmax}^i, L_{\tmin}^i$ with $S_{u,\tmax}, S_{u,\tmin}$, respectively, that is we use the same maximum and minimum values for all the features of the snapshot matrix.  Transformation Eq.~\eqref{eq:normalization} is applied also to the validation and testing sets, but considering as maximum and minimum the values computed over the training set. In order to rescale the reconstructed solution to the original values, we apply the inverse transformation of Eq.~\ref{eq:normalization}.
The configuration of the DL-ROM neural network used in the following test cases is the one given below. We choose a 12-layer DFNN equipped with 50 neurons per hidden layer and $n$ neurons in the output layer, where $n$ represents the dimension of the (nonlinear) reduced  {trial} manifold. The architectures of the encoder and decoder functions are instead  reported in Tables \ref{table_transposed_convolutional_layers_encoder} and \ref{table_transposed_convolutional_layers}.  No activation function is applied at the last convolutional layer of the decoder neural network, as usually done in AEs. 
To solve the optimization problem Eqs.~\eqref{eq:minimization_problem}-\eqref{eq:loss_N}, we use the ADAM algorithm \cite{kingma2015adam}, which is a stochastic gradient descent method computing an adaptive approximation of the first and second momentum of the gradients of the loss function. In particular, it computes exponentially weighted moving averages of the gradients and of the squared gradients. We set the starting learning rate to $\eta = 10^{-4}$, the batch size to 40, and  perform cross-validation in order to tune the hyperparameters of the DL-ROM, by splitting the data in training and validation sets with a proportion 8:2. Moreover, we implement an early-stopping  regularization technique  to reduce overfitting \cite{goodfellow2016deep}, stopping the training if the loss does not decrease over 500 epochs. The maximum number of epochs is set to 10000. As nonlinear activation function, we employ the ELU function \cite{clevert2015fast}. The parameters, weights and biases, are initialized through the He uniform initialization \cite{he2015delving}. The interested reader can refer to \cite{fresca2020POD} 
for a detailed  version of the algorithms used for the training/testing phases. The latter have been carried out on a Tesla V100 32GB GPU by means of the Tensorlow DL framework \cite{abadi2016tensorflow}.

The reduction process requires training the $f_{\phi-\omega}^{DF}$ DFNN as well. The architecture of the NN is the same for both  resonator and  micromirror.
The DFNN receives as input cosine and sine values of a phase $\phi$ and a load multiplier $\beta$ and gives as result a frequency $\omega$.
 Since the solutions are computed with HB method the calculation of the phase response $\phi$ on the first harmonic component is straightforward since it corresponds to $\text{tan}(a_1/b_1)$ with $a_1$ amplitude of the first harmonic sine and $b_1$ amplitude of the first harmonic cosine.
In the next sections, we will plot the results of the $f_{\phi-\omega}^{DF}$ DFNN using the phase $\phi$ instead of the cosine and sine for sake of clarity. 
During the online stage this feature extraction preprocess is embedded in the neural network input, thus only a value for the phase $\phi$ and a forcing level $\beta$ are provided.
The function $f_{\phi-\omega}^{DF} : \mathbb{R}^3 \rightarrow \mathbb{R}$ is a 10-layer DFNN equipped with 64 neurons per hidden layer and 1 neuron in the output layer. We employ the tanh activation function, Adam optimization algorithm, mean squared error metric and batch size equal to 20. The starting learning rate to $\eta = 10^{-3}$ and  perform cross-validation in order to tune the hyperparameters NN, by splitting the data in training and validation sets with a proportion 1:1. 
The parameters, weights and biases, are initialized through the Glorot Normal initialization \cite{glorot2010understanding}.
The early-stopping criterion is used to determine the best fit with a patience coefficient of 50 epochs over a total number of epochs equal to 2000.
The implementation was carried out in  the Tensorlow DL framework \cite{abadi2016tensorflow}.

\begin{figure}[h!]
	\centering
	\includegraphics[width = .8\linewidth]{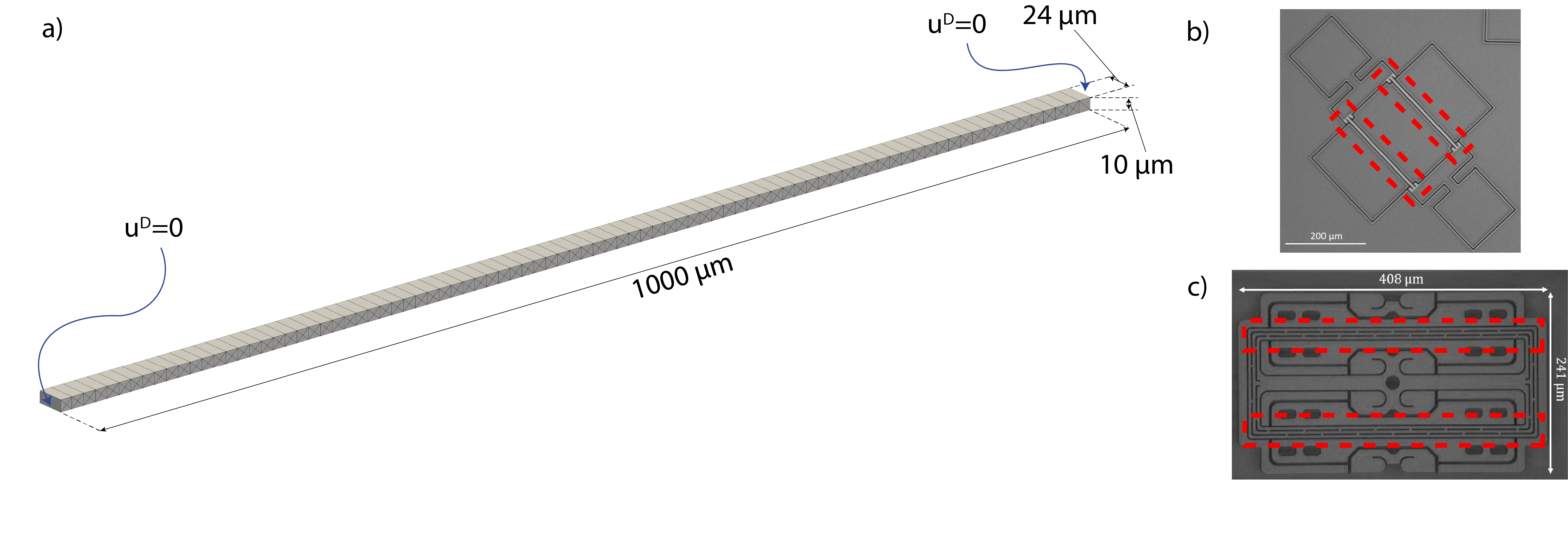}
	\caption{Clamped-clamped beam. Figure a): geometry and mesh of the doubly clamped beam, 3D view. This academic example can be considered as a simplified version of real MEMS resonators like the ones in Figure b) and c) analysed by Zega et al.\cite{jmems20reso}. The red boxes frame the beam resonators in the picture.}
	\label{fig:beam}
\end{figure}

\begin{figure}[h!]
	\centering
	\includegraphics[width = .95\linewidth]{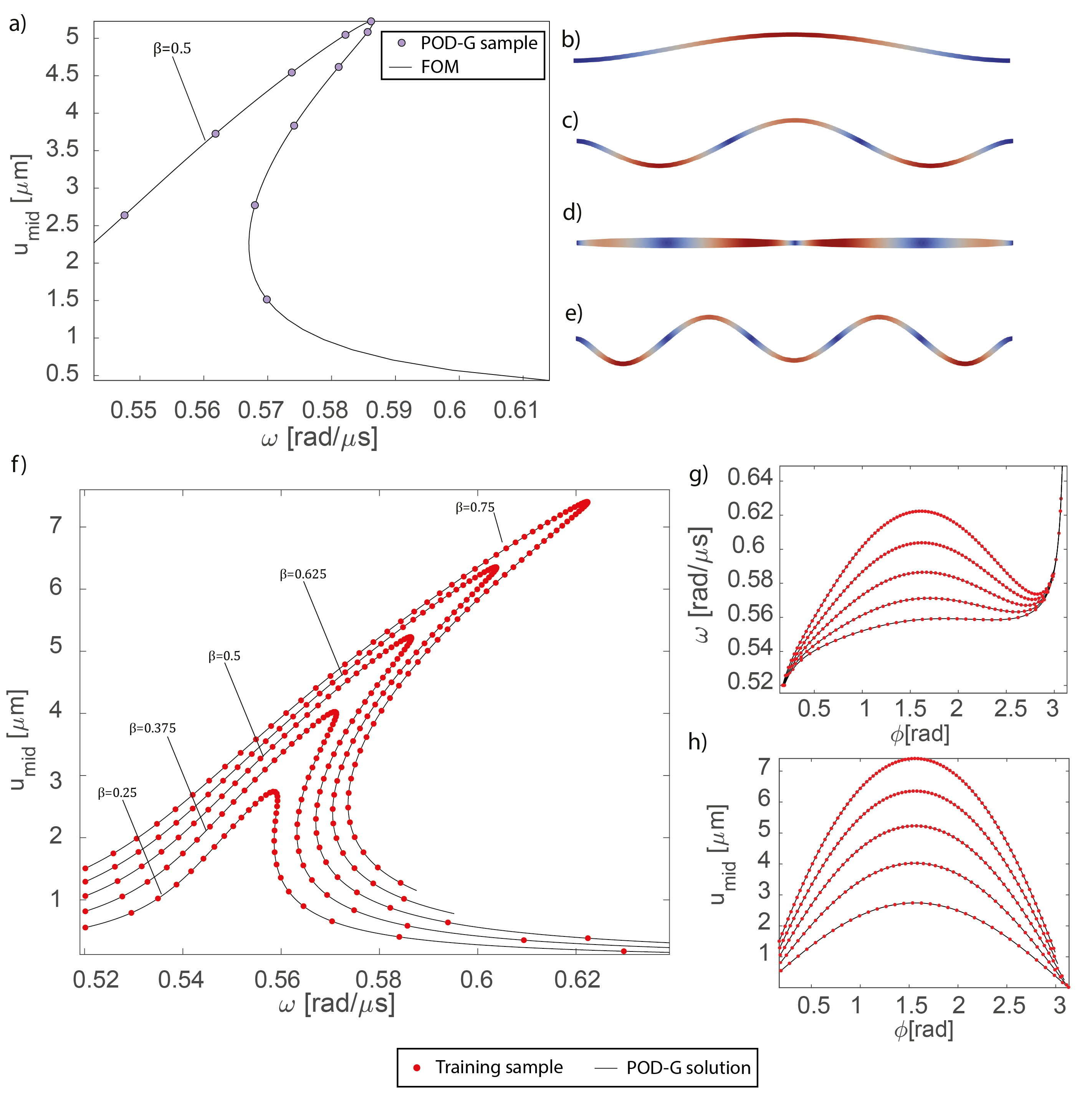}
	\caption{Clamped-clamped beam. Figure a): FOM sampling points ($\beta=0.5$) used to generate the snapshots for the POD-G ROM. Figures b)-e): first four most energetic POD modes obtained from SVD. Figure f): FRFs $u_{\text{mid}}(\omega)$ computed with the POD-G ROM for different load multipliers $\beta$; red markers identify the training dataset used for the DL-ROM (some points have been clipped at higher frequencies). Figures g) and h): plots of the 
	$\omega(\phi)$ and $u_{\text{mid}}(\phi)$ 
	functions corresponding to the FRFs of Figure f), with the same $\beta$ values. }
	\label{fig:beam_POM}
\end{figure}

\begin{figure}[h!]
	\centering
	\includegraphics[width = .95\linewidth]{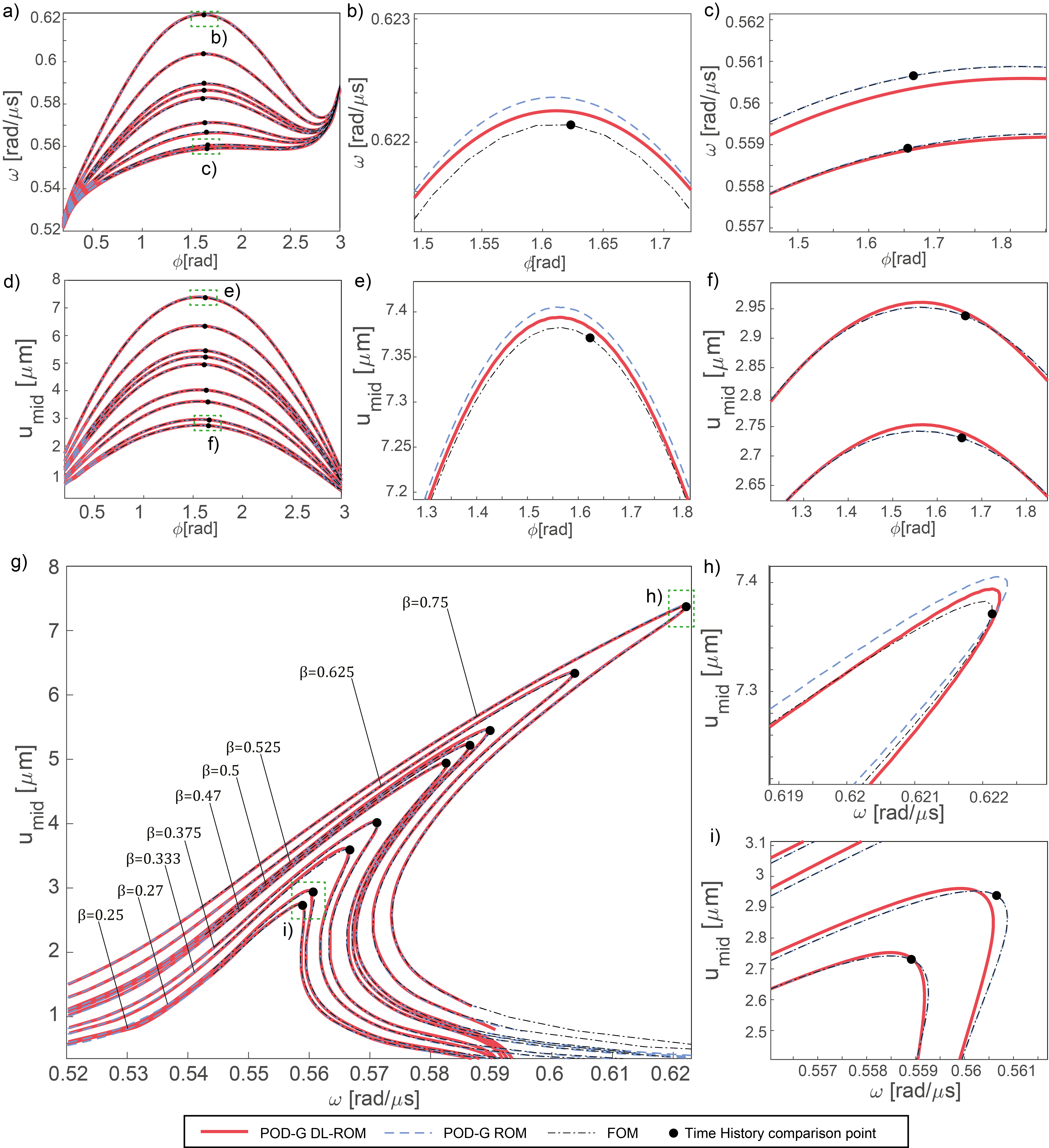}
	\caption{Clamped-clamped beam. Results of the POD-G DL-ROM (red lines) compared with the FOM (black dash-dotted lines) and with the POD-G ROM (dashed blue lines).
    Figure g): final $u_{\text{mid}}(\omega)$ FRF curves for different $\beta$ values. 
    Figures h) and i): enlarged views of the green boxes of figure g). 
    The black dots denote points that will be further investigated in Figure~\ref{fig:time_beam}.
	Figure a): $\omega(\phi)$ curves for the same $\beta$ values of Figure g). 
	Figures c) and d): enlarged views of the green boxes of Figure a).
	Figure d): $u_{\text{mid}}(\phi)$ curves for the same $\beta$ values of Figure g). 
	Figures e) and f): enlarged views of the green boxes of Figure d).} 
	\label{fig:FRF_beam}
\end{figure}

\subsection{Doubly clamped beam resonator}
\label{sec:ccbeam}

Let us consider a doubly clamped beam of length $L=1000$\micr with a rectangular cross-section
of dimensions 10\micr$\times$24\micr, as depicted in Figure~\ref{fig:beam}a). 
Even though this example is proposed primarily to discuss the performance of the proposed technique in the simplest possible setting, nevertheless it also emulates realistic MEMS resonators like those analysed
in \cite{jmems20reso} and depicted in Figures.~\ref{fig:beam}b) and \ref{fig:beam}c). 
The beam is made of isotropic polysilicon \cite{jmems04}, with density $\rho=2330$\,Kg/m$^3$, Young modulus $E=167$\,GPa and  Poisson coefficient $\nu=0.22$.



\begin{figure}[h!]
	\centering
	\includegraphics[width = .8\linewidth]{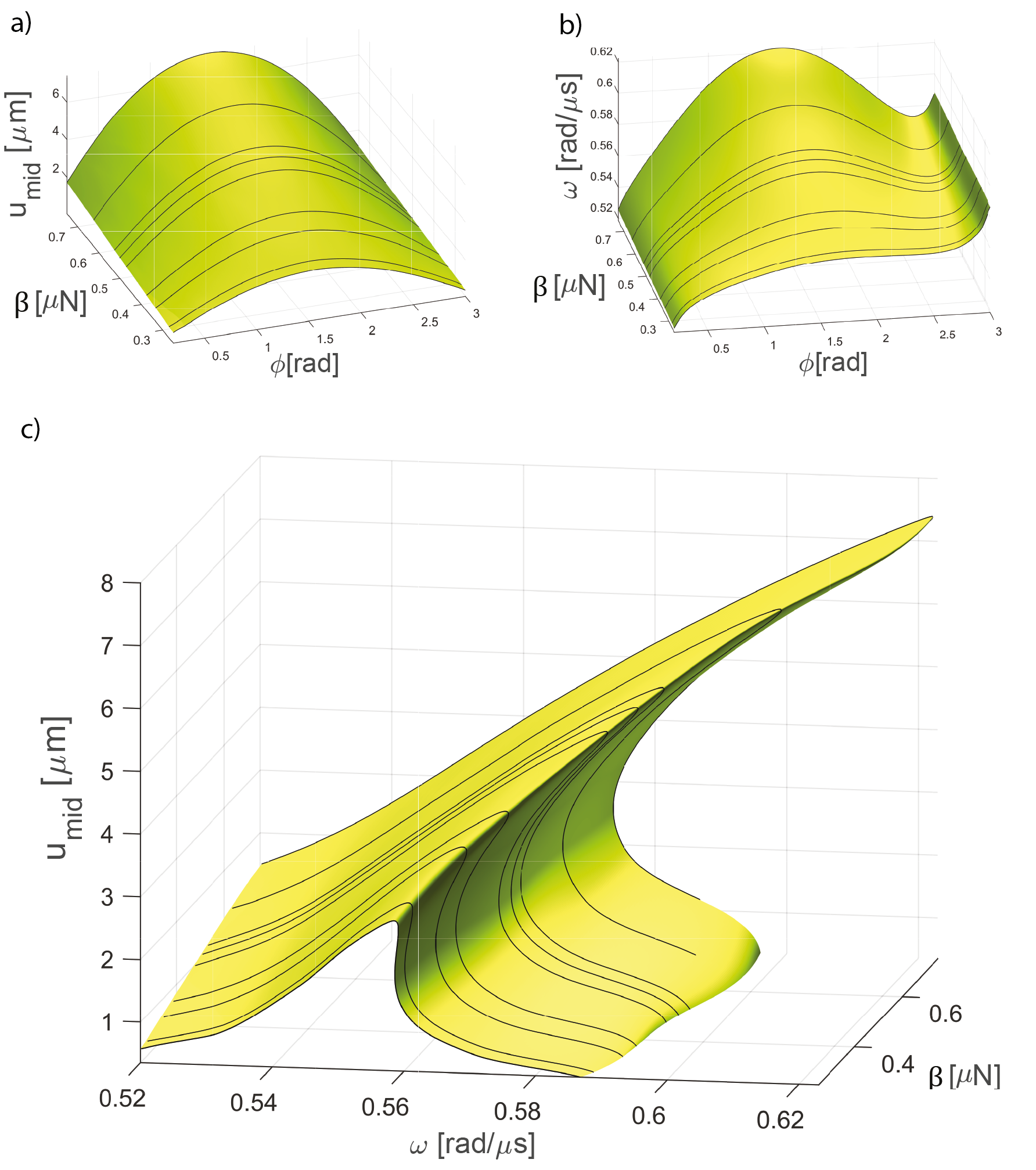}
	\caption{Clamped-clamped beam. Solution envelops generated with the POD-G DL-ROM. 
	Figure a): envelop surface for $u_{\tmid}(\phi,\beta)$. Figure b): envelop surface for $\omega(\phi,\beta)$. Figure c) envelop of the FRFs $u_{\tmid}(\omega,\beta)$. 
	The black lines are the POD-G DL-ROM solutions shown in Figure~\ref{fig:FRF_beam} and are here reported to better highlight the surface curvature.}
	\label{fig:beam_surf}
\end{figure}
\begin{figure}[h!]
	\centering
	\includegraphics[width = .5\textheight]{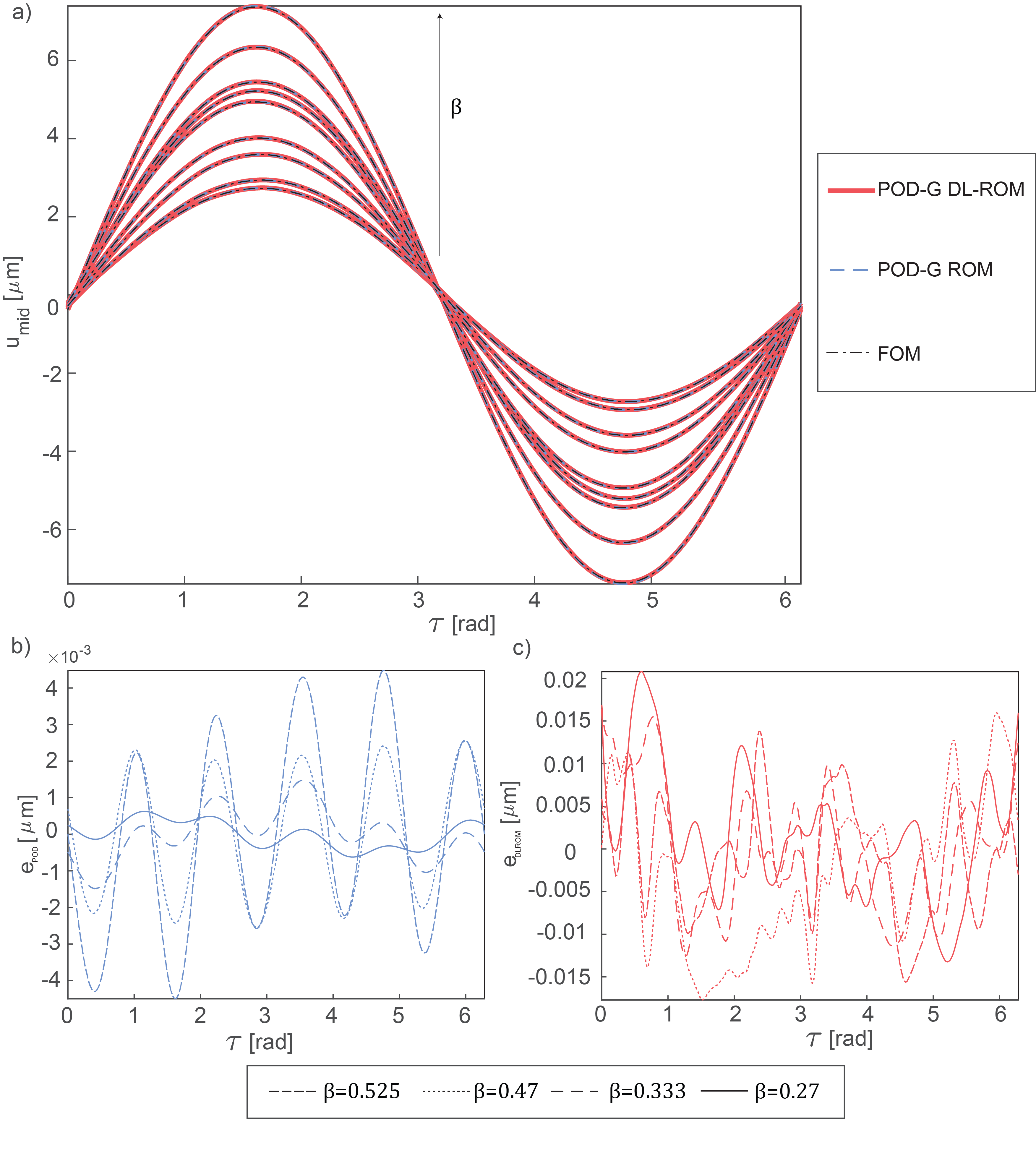}
	\caption{Clamped-clamped beam. Comparison between the periodic responses of the FOM (black line), the POD-G ROM (blue line) and the POD-G DL-ROM (red line) solutions. The sampling points are the ones detailed in Figure~\ref{fig:FRF_beam}. Figure a) depicts the system response along each period. Figures b) and c) represent the error committed by the POD-G ROM and the POD-G DL-ROM with respect to the FOM, computed as the difference between the midspan displacement instantaneous values}
	\label{fig:time_beam}
\end{figure}

\begin{figure}[h!]
	\centering
	\includegraphics[width = .8\linewidth]{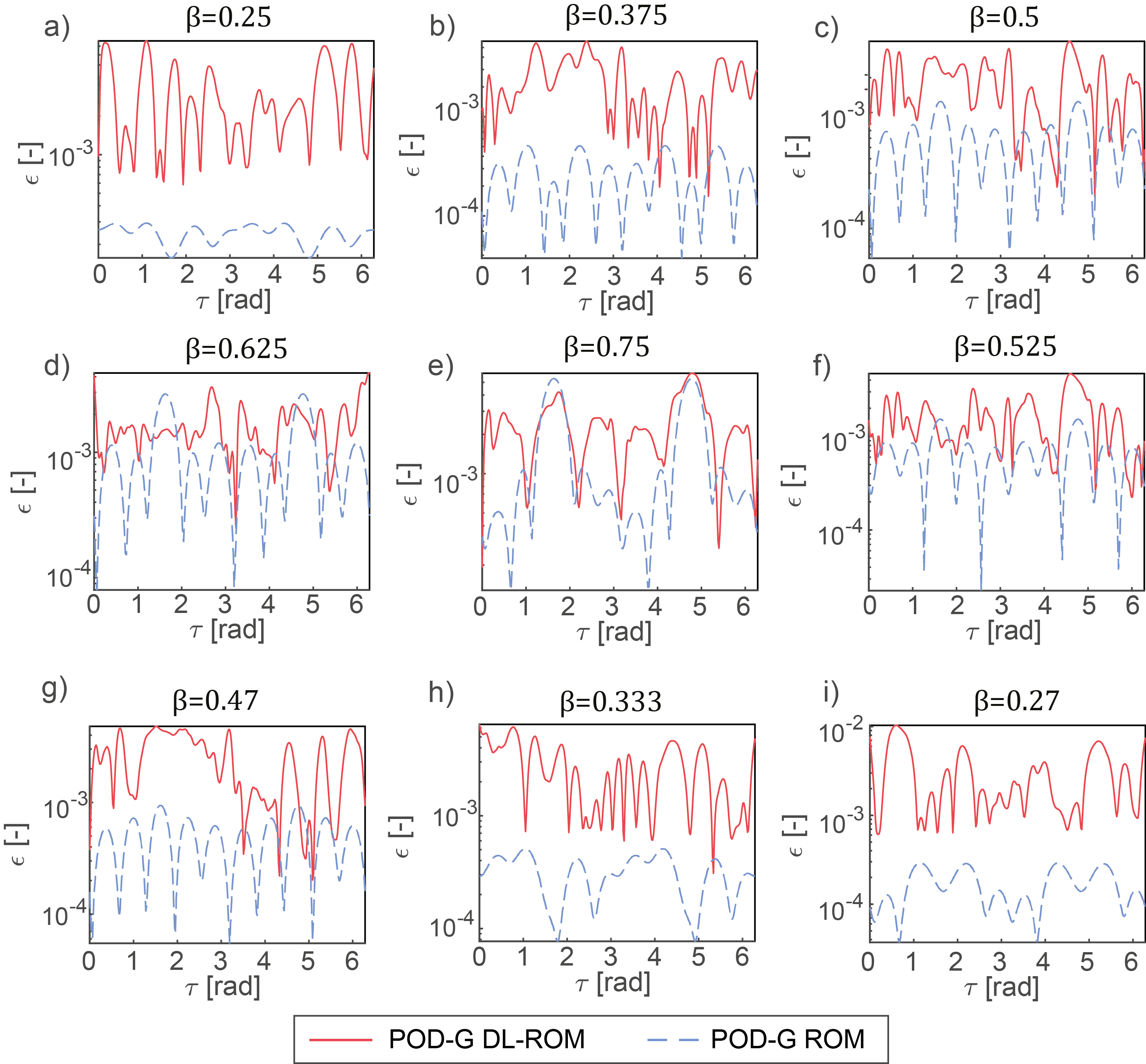}
	\caption{Clamped-clamped beam. Relative error $\epsilon(t)$ (Eq.~\eqref{eq:err}), with respect to the FOM solution, of the POD-G ROM (blue line) and the POD-G DL-ROM (red line) approximations. The sampling points are the ones detailed in Figure~\ref{fig:FRF_beam}. 
	}
	\label{fig:error_beam}
\end{figure}

\begin{figure}[h!]
	\centering
	\includegraphics[width = .5\linewidth]{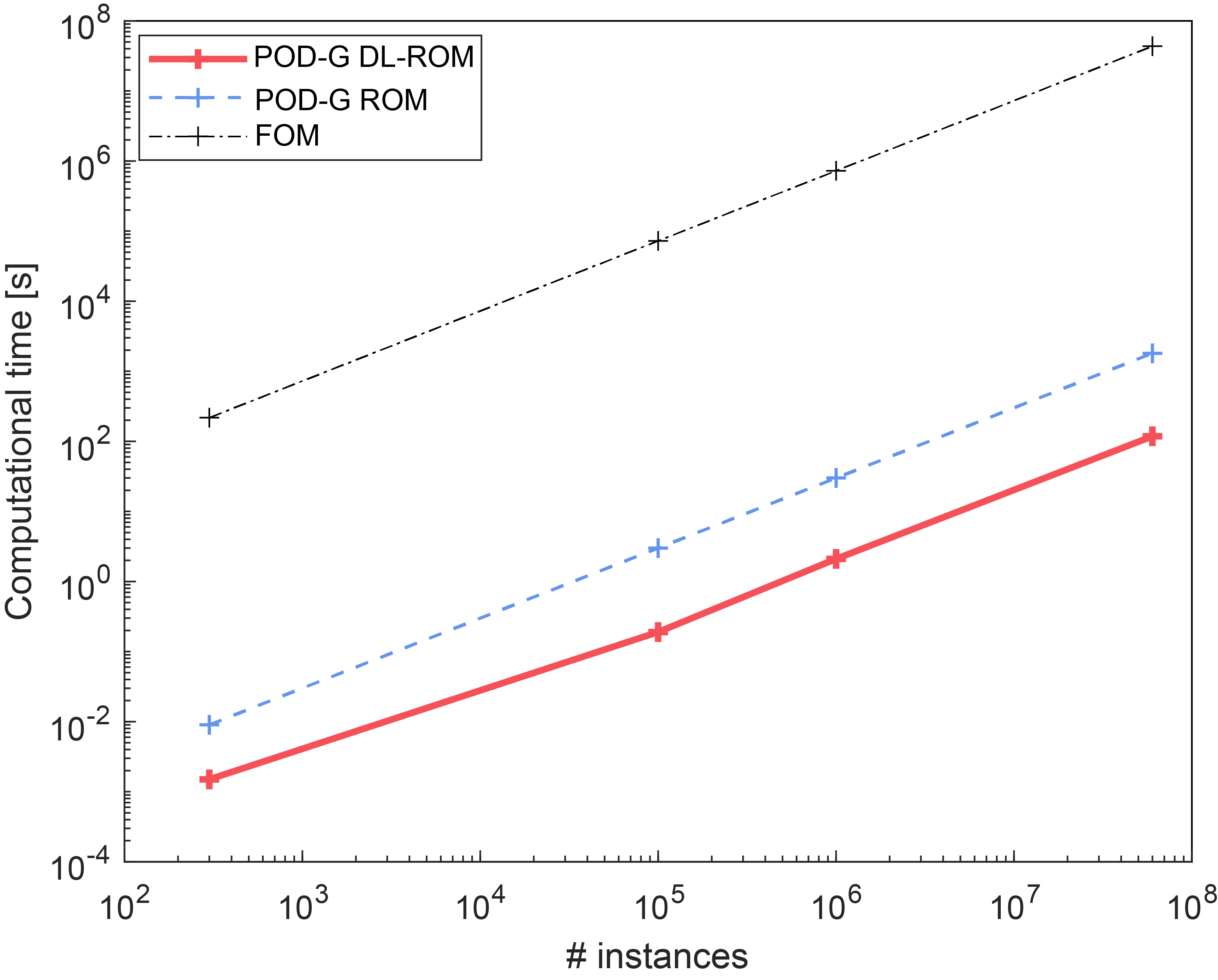}
	\caption{Clamped-clamped beam. Plot of the computational time as a function of the number of instances for the FOM, POD-G ROM and POD-G DL-ROM. All the methods display a linear trend.}
	\label{fig:beam_time}
\end{figure}

Dirichlet boundary conditions ($\bfu^D=\zerou$) are applied on the two opposite sides of the beam. Neumann boundary conditions ($\bff^D=\zerou$) are enforced elsewhere.
The quality factor is assumed equal to $Q=50$. 
The device is excited, for simplicity, by a body load $\bfF=\bfM \bfphi_1 \beta \cos(\omega t)$ 
proportional to the first eigenmode $\bfphi_1$,  with $\bfM$ mass matrix and $\beta$ load multiplier.
The discussion aims at demonstrating the accuracy of the POD-G DL-ROM technique by proposing comparisons with the HB-FOM. 
Hence, in order to keep the FOM computational time at a reasonable level, a rather coarse mesh with 2607 nodes has been employed.
The first five eigenfrequencies are reported in Table~\ref{tab:CC_freq}.

\begin{table}[h!]
	\centering
	\begin{tabular}{ |c|c|c|c|c|c| } 
		\hline
		Eigenmode & 1 & 2 & 3 & 4 & 5\\ 
		\hline
		Frequency  [kHz]& 87.141 & 208.45 & 240.03 & 470.10 & 572.18 \\ 
		\hline
	\end{tabular}
	\caption{Doubly clamped beam: first 5 eigenfrequencies of the FEM model}
	\label{tab:CC_freq}
\end{table}

The POD-G ROM 
basis vectors are  computed using 500 high-fidelity snapshots generated from one period of 10 frequency samples, picked along the FRF corresponding to $\beta=0.5\,\mu$N. The sampling values are represented by the violet circle markers in Figures~\ref{fig:beam_POM}a).
The POD-G ROM is built retaining in the linear trial space only the first 4 most energetic POD modes, depicted in Figures~\ref{fig:beam_POM}b), \ref{fig:beam_POM}c), \ref{fig:beam_POM}d) and \ref{fig:beam_POM}e), respectively.
This number of POD modes represents the minimum required to achieve a good accuracy while keeping the model relatively small \cite{POD21}.

Next, the POD-G ROM is used to generate the training data for the DL-ROM.
The training dataset is given by 309 combinations of phase and forcing. Their positions on the FRFs are indicated by the red markers, in Figure~\ref{fig:beam_POM}f).
The parameter space spanned is $\beta=\{0.25,0.375,0.5,0.625,0.75\}\mu$N$\times \phi=[0.1714:3.1344]$rad. 

To better understand the efficacy of the procedure we compare, 
in Figure~\ref{fig:FRF_beam}, the solutions obtained with the FOM, the POD-G ROM and the POD-G DL-ROM approaches in terms of the plots $\omega(\phi)$, $u_{\tmid}(\phi)$ (with $u_{\tmid}$ maximum midspan displacement) and $u_{\tmid}(\omega)$ (the FRFs) for different $\beta$ values.

The plots of the frequency-phase functions are reported in Figure~\ref{fig:FRF_beam}a)
while enlarged views of specific regions are proposed in  Figures~\ref{fig:FRF_beam}b) and \ref{fig:FRF_beam}c). 
The FOM, POD-G ROM and POD-G DL-ROM solutions are nearly superimposed along the whole curves, some differences being observed where the phase approaches the vertical asymptotes. 
We notice that the POD-G ROM is accurate everywhere with respect to the FOM. 
Only minor differences can be observed at higher amplitudes which can be ascribed to the limited number (4) of POD modes retained in the ROM. 
The POD-G DL-ROM induces a slight increment in the error with respect to the FOM. Nevertheless, its accuracy is still excellent.
The phase-amplitude functions are presented in Figure~\ref{fig:FRF_beam}d), while enlarged views are reported in Figures~\ref{fig:FRF_beam}e) and \ref{fig:FRF_beam}f).
Also in this case FOM, POD-G ROM and POD-G DL-ROM solutions are nearly superimposed, even at the boundaries of the phase parameter space. 
The FRF given by the combination of the previous functions is shown in Figure~\ref{fig:FRF_beam}g), with enlarged views being presented in Figures~\ref{fig:FRF_beam}h) and \ref{fig:FRF_beam}i).
We notice that the resulting FRFs are almost identical and deviate only very close to the limits of the  $\omega(\phi)$ approximations (e.g. close to $\omega=0.59$\,rad/$\mu s$ on the lower branch of the FRF).

All the periodic solutions computed with the POD-G DL-ROM are represented as an envelope of all the FRFs for different values of the forcing in Figure~\ref{fig:beam_surf}.
Such a plot, that would be demanding even with a POD-G ROM, can be completed with the 
POD-G DL-ROM in few seconds and opens the way to the use of this technique in the real time
simulation and optimisation of complex systems.

A further assessment of the POD-G DL-ROM performance concerns time histories along one period.
In Figure~\ref{fig:time_beam}a) we plot the periodic responses of the midspan displacement corresponding to the parameter instances highlighted with the black circle markers in Figure~\ref{fig:FRF_beam}, while Figures~\ref{fig:time_beam}b) and \ref{fig:time_beam}c)
illustrate the error of the POD-G ROM and the POD-G DL-ROM with respect to the FOM.
It is worth stressing that
the procedure reconstructs, for each parameter query and along the whole period, an approximation of the FOM solution everywhere on the structure and not only at selected points.
To quantify the error globally introduced in the spatial domain, we define the instantaneous relative measure:
\begin{equation}
\label{eq:err}
    \epsilon(t)=\frac{||\bfu_h(t)-\tilde{\bfu}_h(t)||}{\sqrt{\frac{1}{N_t}\sum_{k=1}^{N_t}||\bfu_h(t_k)||^2}}
\end{equation}
where $N_t$ denotes the number of time instances considered, $\bfu_h(t)$ is the high-fidelity solution vector at time $t$, $\tilde{\bfu}_h(t)$ is the POD-G DL-ROM reconstruction of the FOM solution and $|| \cdot ||$ denotes the L2-norm at a given time $t$. The denominator of Eq.~\eqref{eq:err} provides a reference value for the normalization given by the average of the L2-norm of the solution along the whole period.
The relative errors along each period, always close to or lower than $\epsilon=10^{-3}$,  
are reported in Figure~\ref{fig:error_beam}.

\begin{table}[h!]
	\centering
	\begin{tabular}{ |c|c|c|c|c|c| } 
		\hline
		number of instances & $T_{\tFOM}$ & $T_{\tPOD}$ & $T_{\tDL}$  & $T_{\tFOM}/T_{\tDL}$ & $T_{\tPOD}/T_{\tDL}$\\ 
		\hline
		$300$ & $212 s$ & $0.009 s$ &  $0.0015 s$  & $1.41\times 10^{5}$ & 6.0 \\
		\hline
		$10^5$ & $20.13 h$ & $3.2 s$ & $0.18 s$  & $4.0\times 10^{5}$ & 17.7\\
		\hline
		$10^6$ & $8.4 d$ & $31.5 s$ & $2.0 s$  & $3.62\times 10^{5}$ & 15.7\\
		\hline
		$6\times10^7$ & $503 d$ & $1890 s$ &  $115 s$  & $3.77\times 10^{5}$ & 16.4\\
		\hline
	\end{tabular}
	\caption{Clamped-clamped beam. Computational cost of the FOM
	($T_{\tFOM}$), POD-G ROM ($T_{\tPOD}$) and POD-G DL-ROM ($T_{\tDL}$)
	and speedup of the reduction techniques.}
	\label{tab:beam_time}
\end{table}

Next, we compare the computational time required to solve the FOM with the POD-G ROM and the POD-G DL-ROM testing times, i.e. the time needed to compute the ROM solution for a given number of parameter instances. The results are summarised in 
Table~\ref{tab:beam_time} and the computational cost trend is depicted in Figure~\ref{fig:beam_time}.

We notice that the $T_{\tFOM}$ explodes when increasing  the number of instances while 
the POD-G ROM ($T_{\tPOD}$) can drastically reduce the computing time 
making it a suitable approach when there is no need of truly real-time applications.
As extensively discussed in \cite{POD21}, 
in this case a robust ROM can be build with only 4 POD modes.
Nevertheless, even in this extremely simple example,  
the POD-G DL-ROM  testing time ($T_{\tDL}$) highlights the impressive performance of the technique, which guarantees a 
reduced effort, with respect to the FOM, of a factor $10^5$.

\subsection{MEMS Micromirror}
\label{sec:micromirror}

Micromirrors represent one of the most promising examples of the next generation of MEMS
which will deeply impact our lives and enable a sustainable evolution of the Internet of Things. 
Among their virtually countless applications, one can cite pico-projectors for Augmented Reality (AR) lenses \cite{ARlenses} and 3D scanners for Light Detection and Ranging (LiDAR) in autonomous driving.

\begin{figure}[h!]
	\centering
	\includegraphics[width = .95\linewidth]{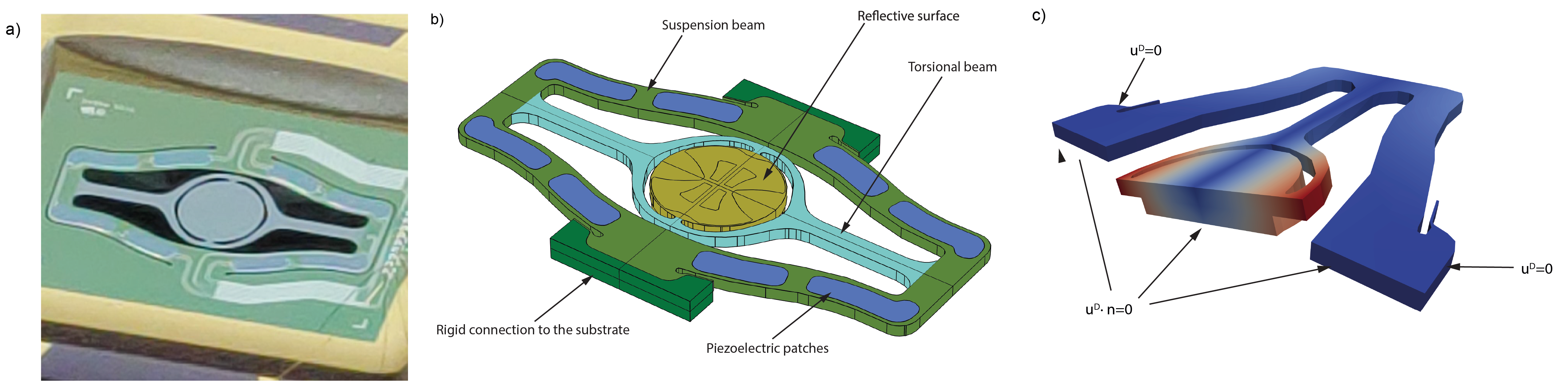}
	\caption{Micromirror. Figure a): optical picture of the micromirror. 
		Figure b): Schematic view of the layout with few details on the device components \cite{actuators21,jmemsmirror}. Figure c): 
		third (torsional) eigenmode that is actuated during operations and the boundary conditions used in the simulation of half of the device.}
	\label{fig:perseus}
\end{figure}

\begin{figure}[h!]
	\centering
	\includegraphics[width = .95\linewidth]{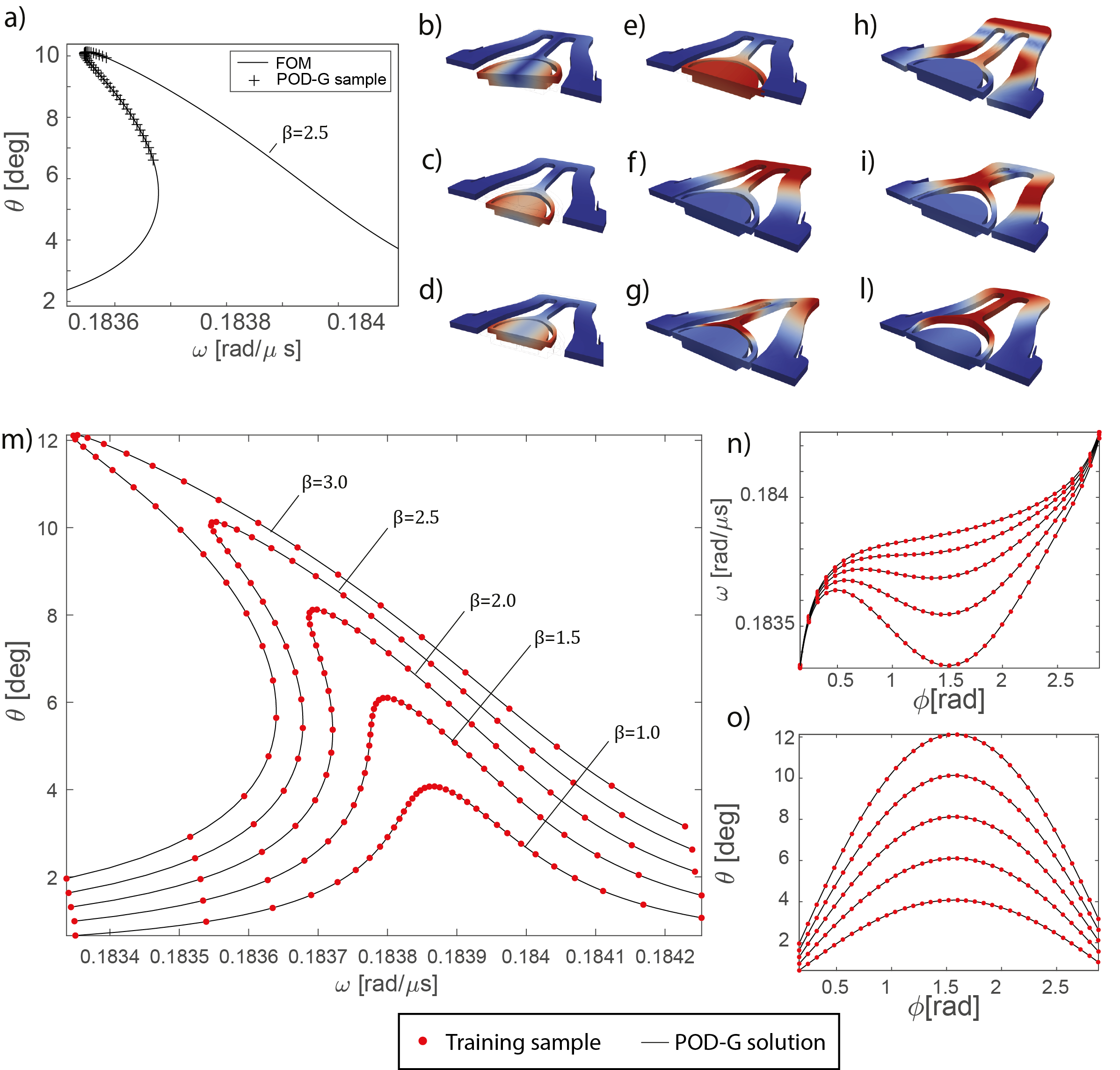}
	\caption{Micromirror.  Figure a): sampling points on the FOM FRF used to build the POD-G ROM subspace bases. Figure b)-l): 9 most energetic POD modes used in the reduced basis. 
	Figure m): FRFs $\theta(\omega)$ computed with the POD-G ROM for different $\beta$ values. 
	The training sampling points used to generate the snapshots for the DL-ROM construction are highlighted with red markers. Some points at higher frequencies are clipped. Figure n) and o):  $\omega(\phi)$ and $\theta(\phi)$ functions corresponding to the FRFs of Figure m), with the same $\beta$ values.}
	\label{fig:perseus_POM}
\end{figure}

\begin{figure}[h!]
	\centering
	\includegraphics[width = .95\linewidth]{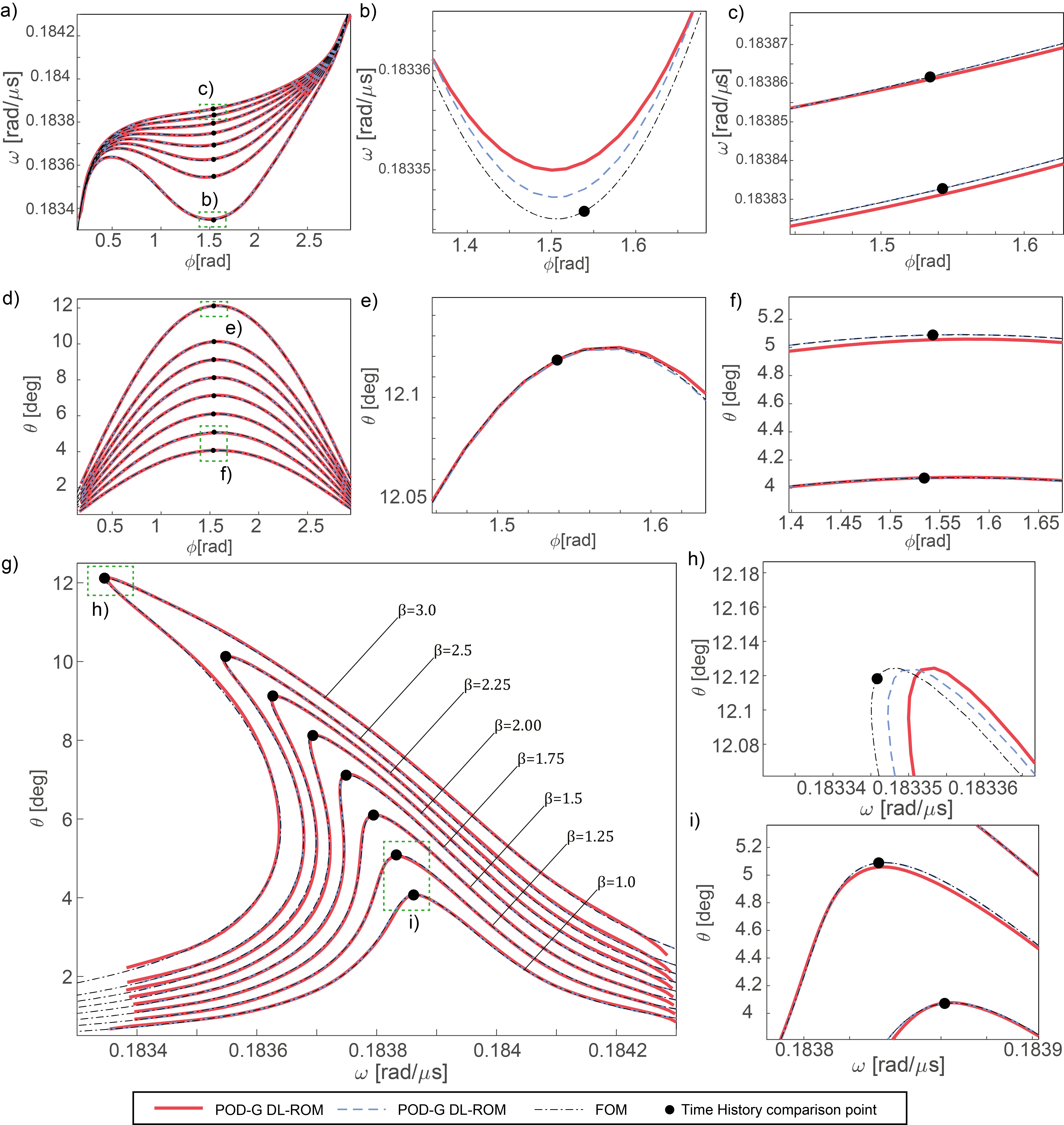}
	\caption{Micromirror: Results of the POD-G DL-ROM (red lines) compared with the FOM (black dash-dotted lines) and with the POD-G ROM (dashed blue lines). Figures a) and d): frequency $\omega(\phi)$ and angle $\theta(\phi)$ plots, respectively, for various forcing levels. 
    Figures b) and c) are enlarged views of $\omega(\phi)$ in Figure a).
	Figures e) and f) are enlarged views of $\theta(\phi)$ in Figure b).}
	\label{fig:FRF_perseus}
\end{figure}

\begin{figure}[h!]
	\centering
	\includegraphics[width = .8\linewidth]{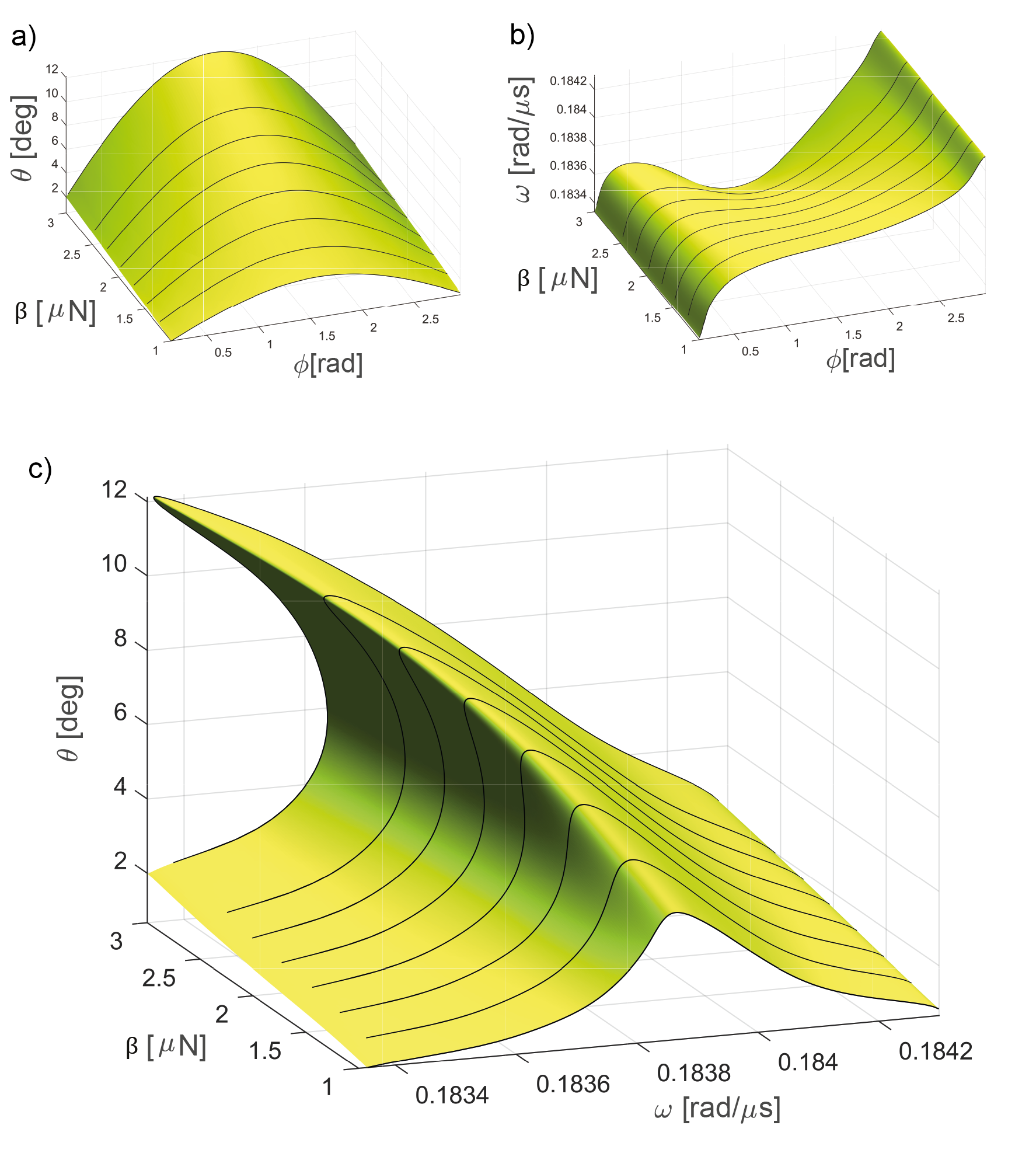}
	\caption{Micromirror. Envelop of the solutions reconstructed with the POD-G DL-ROM. 
	Figure a): $\theta(\phi;\beta)$ envelop. Figure b): $\omega(\phi;\beta)$ envelop. 
	Figure c): $\theta(\omega;\beta)$ envelop. The black line are the POD-G DL-ROM solutions shown in Figure~\ref{fig:FRF_beam} and are here reported to better represent the surface curvature.}
	\label{fig:perseus_surf}
\end{figure}

Because of the inertial and geometrical effects triggered by large rotations, micromirrors are intrinsically nonlinear and the prediction of their dynamic behaviour is essential to guarantee a proper design and control of the mirror during the online stage.
Only recently their numerical simulation has been tackled by Opreni et al.\cite{actuators21}
with a dedicated large scale FOM based on HB method.
It is worth stressing that the generation of a suitable ROM 
is a tough challenge even for the most advanced and recent techniques, like the Direct Normal Form approach,
mainly because the torsional mode is not the lowest-frequency
one (it is the third) and is not well separated from the other modes. 
Indeed, the quadratic formulation of the Direct Normal Form approach 
implemented by Opreni et al. \cite{NNM21}  fails and a high order expansion is required,
as remarked by Vizzaccaro et al.\cite{vizzaccaro2021high}.
Similar difficulties have been experienced with the Implicit Condensation approach 
\cite{ijnm19}. 

The mirror considered, fabricated by STMicroelectronics, is illustrated in Figure~\ref{fig:perseus}. The mirror plate is  suspended to a gimbal connected with a torsional beam along the rotation axis and two suspension beams on each side.
The mirror is assumed to be made of isotropic polysilicon \cite{jmems04}, with density $\rho=2330$\,Kg/m$^3$, Young modulus $E=167$\,GPa and  Poisson coefficient $\nu=0.22$. Since the central plate is very stiff, we adopt its angle of rotation $\theta$ as reference for the FRFs amplitude and for the time response.
Thanks to symmetry, only half of the micromirror is modelled with the FEM and a total of 9732 dofs. The Dirichlet boundary conditions are imposed on the substrate (the dark green areas in
Figure~\ref{fig:perseus}b) and on the symmetry plane. On the remaining boundaries, zero traction Neumann boundary conditions are imposed.
The first five eigenfrequencies obtained from a linear FEM eigenvalue analysis are listed in Table~\ref{tab:mirr2_freq}; the torsional mode of interest herein, the third one, has a frequency of $29 271$\,Hz.
A quality factor equal to $Q=1000$ is considered.

\begin{table}[h]
	\centering
	\begin{tabular}{ |c|c|c|c|c|c| } 
		\hline
		Eigenmode & 1 & 2 & 3 & 4 & 5\\ 
		\hline
		Frequency  [kHz]& 11.080 & 18.533 & 29.271 & 41.667 & 68.848 \\ 
		\hline
	\end{tabular}
	\caption{Micromirror: first five eigenfrequencies computed with a linear FEM}
	\label{tab:mirr2_freq}
\end{table}

\begin{figure}[h!]
	\centering
	\includegraphics[width = .45\textheight]{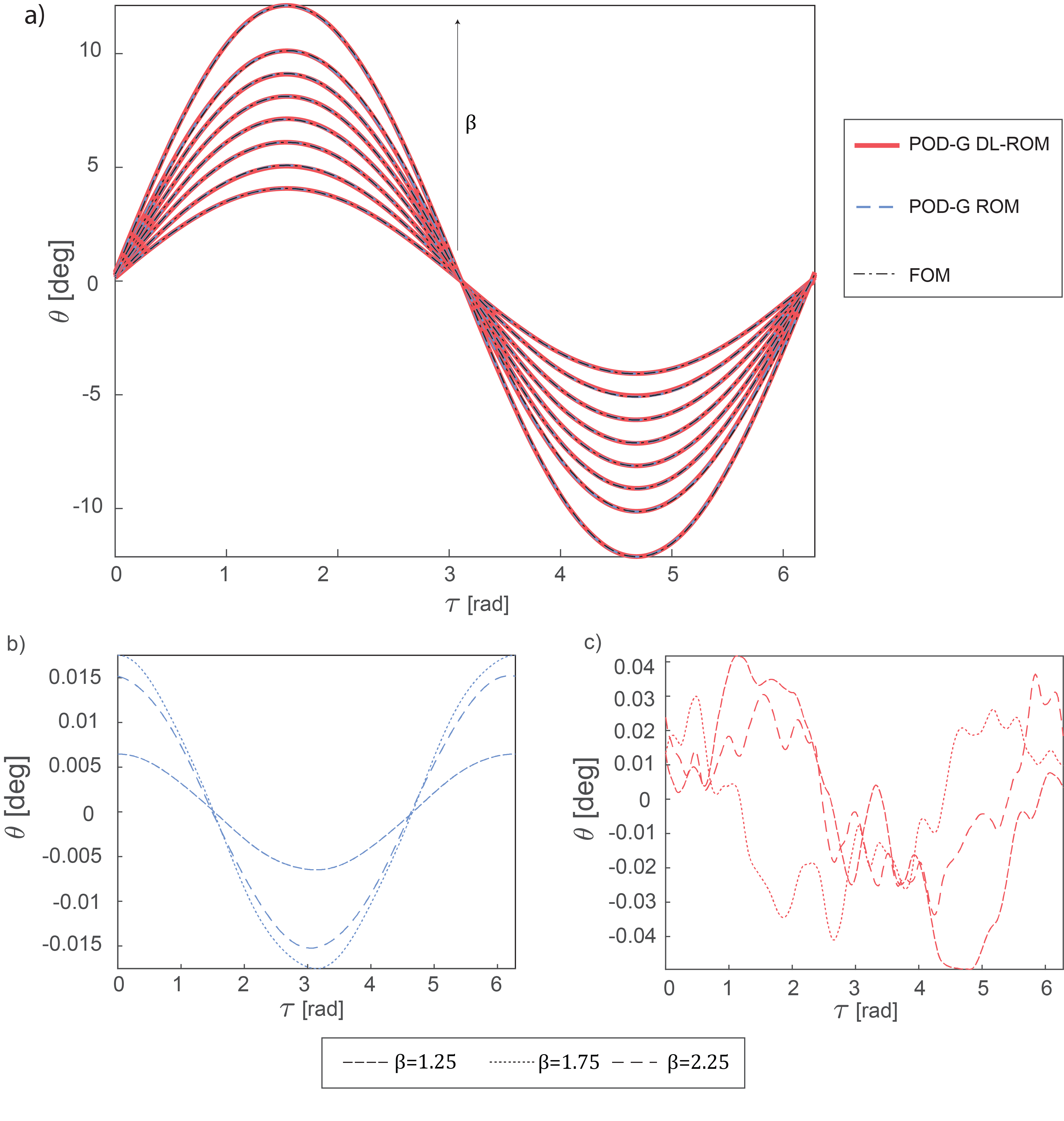}
	\caption{Micromirror. Comparison between the periodic responses of the FOM (black line), the POD-G ROM (blue line) and the POD-G DL-ROM (red line) solutions. The sampling points are indicated by black dots in Figure~\ref{fig:FRF_perseus}. Figure a): evolution of the rotation angle over one period. Figures b) and c): angle error of the POD-G ROM and the POD-G DL-ROM solutions with respect to the FOM.}
	\label{fig:time_perseus}
\end{figure}

\begin{figure}[h!]
	\centering
	\includegraphics[width = .6\linewidth]{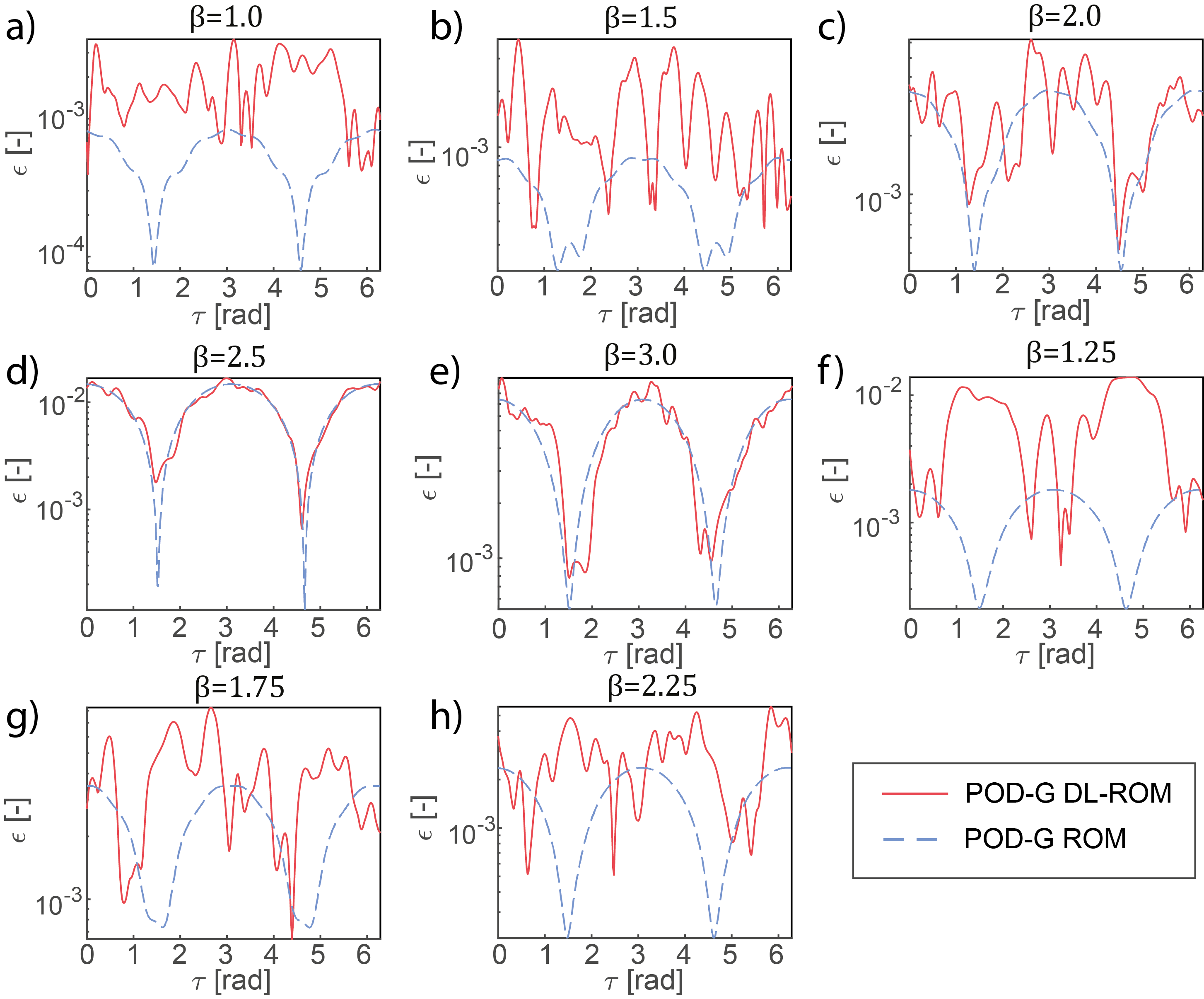}
	\caption{Micromirror. Relative error $\epsilon(t)$ (Eq.~\eqref{eq:err}) with respect to the FOM solution for the POD-G ROM (blue line) and the POD-G DL-ROM (red line) approximations. The sampling points are the ones detailed in Figure~\ref{fig:FRF_perseus}. 
	}
	\label{fig:err_perseus}
\end{figure}

\begin{figure}[h!]
	\centering
	\includegraphics[width = 1.\linewidth]{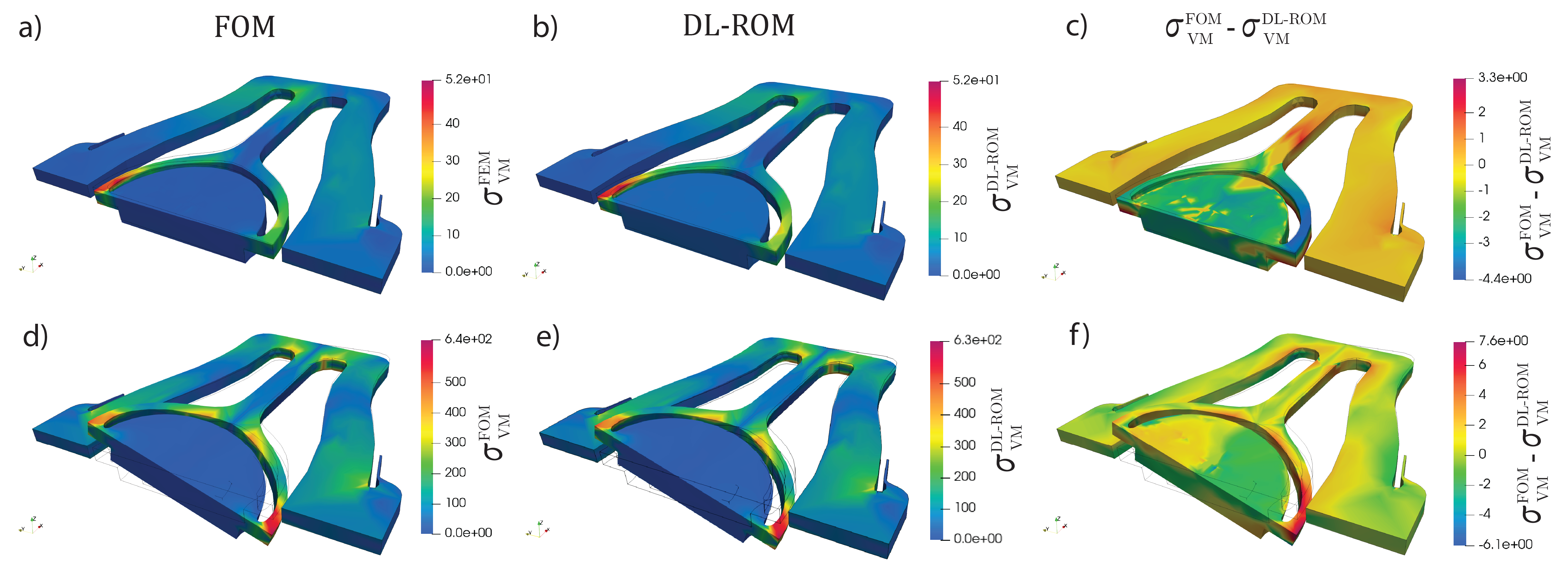}
	\caption{Contour plot of the Von Mises stress $\sigma_{\text{VM}}$  at the peak of the FRF for $\beta=2.25$ and two different instants of the periodic response i.e.\ $\tau=\pi/2$ and $\tau=\pi$, with $\tau$ 
	normalised time. Figures a) and d): field computed with the FOM. Figures b) and e): field computed with the POD-G DL-ROM. 
	Figures c) and f): difference between the two reconstructions. The displacement of the mirror has been amplified by a factor 5}	
	\label{fig:sforzo_perseus}
\end{figure}

The micromirror is inspired by the one analysed by Opreni et al.\cite{actuators21}, nevertheless we simplify the forcing mechanism by replacing the piezoelectric actuation with a body load proportional to the third torsional eigenmode  $\bfF=\bfM \bfphi_3 \beta \cos(\omega t)$. Such a simplification allows focusing only on the application of the POD-G DL-ROM technique to the modelling of geometric and inertia nonlinearities.

As for the previous example, an HB-FOM is used to generate the high-fidelity solutions to be processed with SVD.
We collect 2000 snapshots obtained for $\beta=2.5\,\mu$N and 40 frequency values distributed as illustrated by the cross markers in Figure~\ref{fig:perseus_POM}a).
To achieve a good balance between subspace dimension and quality of the solution, the 
9 POD modes illustrated in Figures~\ref{fig:perseus_POM}b)-\ref{fig:perseus_POM}l) are retained in the reduced basis.
A critical analysis of the quality of the POD solution in terms of the number of POD modes utilised can be found in the work from Gobat et al.\cite{POD21}.

The POD-G ROM is next used to generate the low-fidelity snapshots for the DL-ROM and the $f_{\phi-\omega}^{DF}$ DFNN.
Both the training datasets contain 175 combinations of phase and forcing, and are denoted in Figure~\ref{fig:perseus_POM}m) by the red markers.
The parameter space spanned is $\beta=\{1.0, 1.5, 2.0, 2.5, 3.0\}{\mu}\text{N}\times \phi=[0.1604:2.9211]\,\text{rad}$.
The trained DL-ROM is finally applied to build the FRF curves for all the desired  forcing levels as plotted in Figure~\ref{fig:FRF_perseus}.
Similarly to what remarked in the previous application to the clamped-clamped beam,
the POD-G DL-ROM reproduces the dynamic behaviour of the system with remarkable accuracy and errors are  limited to regions corresponding to the lower and upper bounds of the phase, of little engineering interest.
Also in this case, the POD-G DL-ROM can be used to produce an almost continuous envelop of the system FRFs, as shown in Figure~\ref{fig:perseus_surf}.
The same accuracy is displayed by the time evolution of the angle, reported in Figure~\ref{fig:time_perseus}, corresponding to the parameter combinations highlighted by black dots in Figure~\ref{fig:FRF_perseus}.
The error of the POD-G ROM and the POD-G DL-ROM with respect to the FOM is depicted in  Figure~\ref{fig:time_perseus}b) and c) respectively.
Finally, the relative error $\veps$, Eq.~\eqref{eq:err} 
is plotted in Figure~\ref{fig:err_perseus} for various load multipliers $\beta$.



Starting from the output $\tilde{\bfu}_N$ of the decoder, the FOM approximate solution can be reconstructed using the POD-G modes as:
$\tilde{\bfu}_h(t;\bfmu) = \mathbf{V} \tilde{\bfu}_N(t;\bfmu)$.
This field can be later manipulated with standard FEM tools to generate 
any quantity of interest, like e.g.\ stresses which are often required to monitor strength criteria. 
As an example,  the Von Mises stress $\sigma_{\text{VM}}$ 
is plotted in  Figure~\ref{fig:sforzo_perseus} 
where the FOM solution is compared with the prediction of the POD-G DL-ROM at two different time instants in Figure~\ref{fig:err_perseus}. A quantitative inspection shows that an excellent accuracy is preserved. Some artefacts appear in the plots which are however associated with the rather coarse mesh utilised in the example.

\begin{figure}[h]
	\centering
	\includegraphics[width = .5\linewidth]{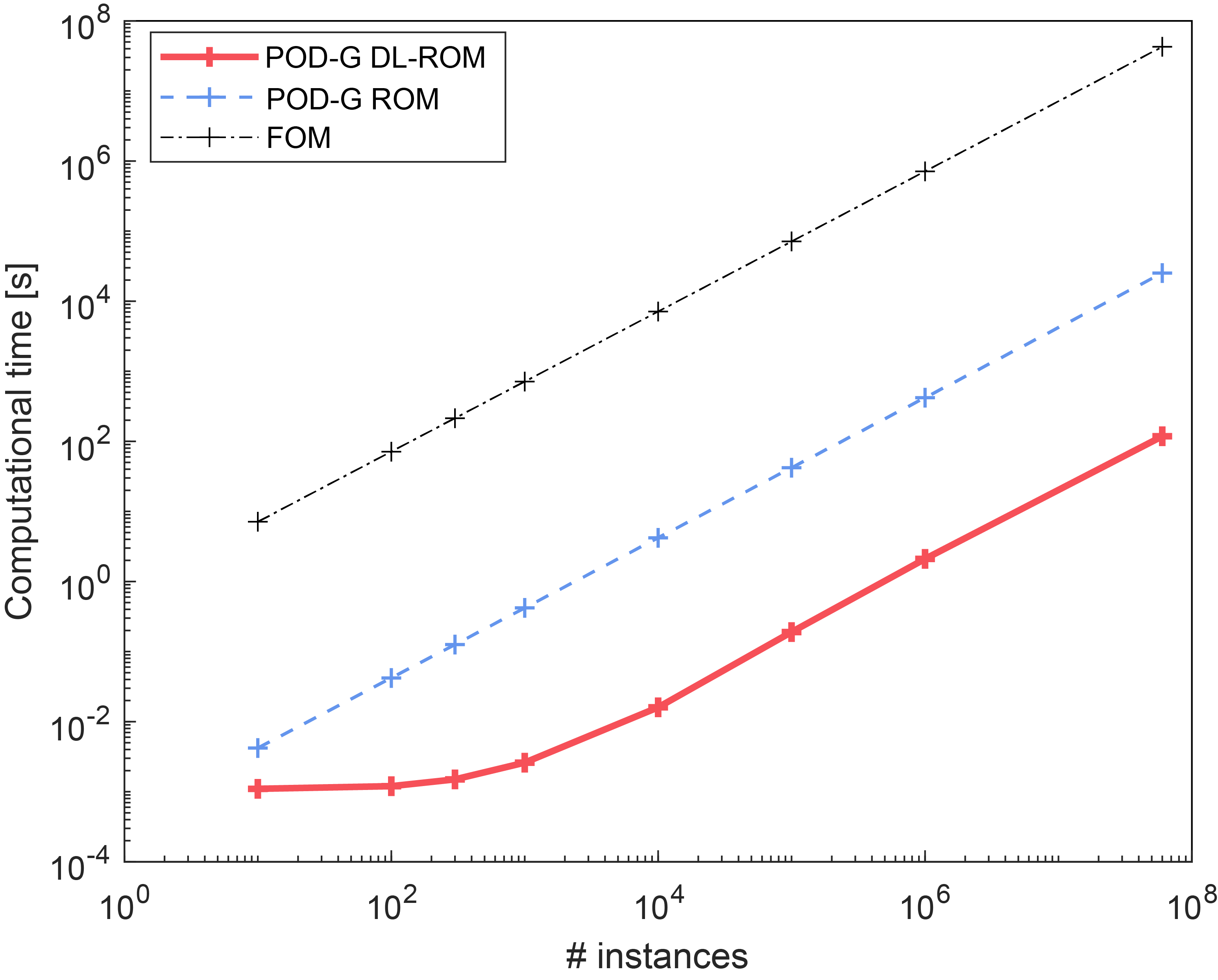}
	\caption{Micromirror. Plot of the time required to compute a given number of instances for the three solution level: FOM, POD-G ROM and POD-G DL-ROM. All the methods display a linear trend. When a few instances are considered the POD-G DL-ROM displays a plateu before the linear trend. This behaviour is a consequence of the parallelization, under a certain threshold the number of cores overcomes the number of instances thus the time needed to get the result is nearly the same.}
	\label{fig:perseus_time}
\end{figure}

Finally, we compare the computational cost of the FOM, the POD-G ROM and the POD-G DL-ROM solutions. The results are summarized in Table~\ref{tab:mirr_time} and the computational cost trend is depicted in Figure~\ref{fig:perseus_time}.
Contrary to the clamped-clamped beam example, the larger POD basis dimension (9 as opposed to 4) 
yields as expected an increase of the ratio $T_{\tPOD}/T_{\tDL}$ which puts in better evidence the interest of the proposed approach over other model order reduction techniques for real-time applications.

\begin{table}[h]
	\centering
	\begin{tabular}{ |c|c|c|c|c|c| } 
		\hline
		number of instances & $T_{\tFOM}$ & $T_{\tPOD}$ & $T_{\tDL}$  & $T_{\tFOM}/T_{\tDL}$ & $T_{\tPOD}/T_{\tDL}$\\ 
		\hline
			$10^1$ & 7.09 s & 0.0042 &0.0011 s & $6.4\times 10^{3}$ & 3.81\\
		\hline
			$10^2$ & 1.18 m & 0.042 s &0.0012 s & $5.9\times 10^{4}$ & 35\\
		\hline
			$3\times 10^2$ & 3.54 m & 0.12 s & 0.0015 s & $1.41\times 10^{5}$  & 80\\
		\hline
			$10^3$ & 11.82 m & 0.42 s & 0.0026 s & $2.72\times 10^{5}$  & 161\\
		\hline
			$10^4$ & 1.97 h & 4.23 s & 0.016 s  & $4.43\times 10^{5}$  & 264\\
	    \hline
			$10^5$ & 19.7 h & 42.3 s &0.19 s  &  $3.73\times 10^{5}$ & 222\\
		\hline
			$10^6$ & 8.2 d & 7 m &2.1 s  & $3.37\times 10^{5}$ & 200\\
		\hline
		    $6\times 10^7$ & 494 d & 7.05 h & 118 s  & $3.6\times 10^{5}$ & 215\\	
		\hline
	\end{tabular}
	\caption{Micromirror. Computational cost of the FOM
	($T_{\tFOM}$), POD-G ROM ($T_{\tPOD}$) and POD-G DL-ROM ($T_{\tDL}$)
	and speedup of the reduction techniques.}
	\label{tab:mirr_time}
\end{table}

\section{Conclusions}
\label{sec:conclusions}

In this paper we have proposed a non-intrusive Deep Learning-based ROM 
tailored to describe both the trial manifold and the reduced dynamics of complex mechanical systems showing inertia and geometric nonlinearities. 
In order to reduce the computational cost of the procedure, a POD-Galerkin ROM has been first generated
using a limited amount of FOM snapshots and the training of the 
DL-ROM has been performed using cheaper POD-G solutions covering the whole parameter range. 

The encoder function of a convolutional autoencoder (AE) has been used 
to map the system response onto a low-dimensional representation, while the decoder part models the reduced nonlinear trial manifold and allows to reconstruct a posteriori an approximation of the FOM response.
Finally, in order to describe the system dynamics on the manifold, a Deep Forward Neural Network has been utilized and trained at the same time as the AE.

The proposed technique has proved extremely efficient as it models highly nonlinear problems by identifying the manifold underlying the dynamics in a completely data-driven, black-box and non-intrusive way, provided that the sampled data span the parameter space of interest. 

This strategy has been implemented and benchmarked on both an academic example 
and a real industrial application like a MEMS micromirrors showing
softening response and multiplicity of solutions. In both cases a truly real-time
ROM is achieved that preserves 
the capability to retrieve a posteriori, through the decoder part of the AE, an 
estimate of any solution field available in the corresponding FOM.

In both applications, the POD-G DL-ROM has reproduced the dynamic behaviour of the system with remarkable accuracy, inducing limited errors in the regions corresponding to the lower and upper bounds of the phase range.
Concerning the computational cost, we have compared the POD-G DL-ROM both against 
the original FOM solved by HB procedures and against the POD-G ROM. In the former case the gain
is always impressive, and in the latter, as soon as the POD trial-space dimension increases
in realistic examples, improvements of two orders of magnitude have been observed.


\section*{Acknowledgments}
SF, GG and AM have been supported by Fondazione Cariplo, Grant n. 2019-4608.

\bibliographystyle{plain}     

\bibliography{DLROM_2021}

\end{document}